\documentclass[12pt]{amsart}

\usepackage{accents}
\usepackage{appendix}
\usepackage{amsfonts}
\usepackage{amsmath}
\usepackage{amssymb}	
\usepackage{amsthm,bm}
\usepackage{array,booktabs,multirow}
\usepackage{braket}
\usepackage{centernot}
\usepackage{cite}
\usepackage{comment}
\usepackage{dsfont}
\usepackage{mathalfa}
\usepackage[shortlabels]{enumitem}
\usepackage{etoolbox}
\usepackage{float}
\usepackage[hang, flushmargin]{footmisc}
\usepackage{latexsym}
\usepackage{lipsum}
\usepackage{needspace}
\usepackage{tikz}
\usepackage{hyperref}
\usetikzlibrary{matrix,arrows}
\usepackage{diagbox}
\usepackage{adjustbox}

\makeatletter
\def\namedlabel#1#2{\begingroup
   \def\@currentlabel{#2}%
   \label{#1}\endgroup
}
\makeatother

\theoremstyle{plain}
\newtheorem{thm}{Theorem}[section]

\newtheorem{lem}[thm]{Lemma}
\newtheorem{prop}[thm]{Proposition}

\theoremstyle{definition}

\newtheorem{defn}[thm]{Definition}

\theoremstyle{remark}

\setlist[enumerate,1]{leftmargin=2em}

\def\A{\mathbf A}

\def\C{\mathbb C}
\def\End{{\rm End}}

\def\I{\mathcal I}
\def\i{\mathbf i}

\def\N{\mathbb N}

\def\R{\mathbf R}
\def\T{\mathbf T}

\def\Z{\mathbb Z}

\def\ds{\displaystyle}

\def\U{U(\mathfrak{sl}_2)}

\def\Ue{U(\mathfrak{sl}_2)_e}
\def\e{\varepsilon}

\newcommand{\floor}[1]{\left\lfloor #1 \right\rfloor}
\newcommand{\ceil}[1]{\left\lceil #1 \right\rceil}

\title[An algebra homomorphism from $\Re$ into $\U$]{Johnson graphs as slices of a hypercube and an algebra homomorphism from the universal Racah algebra into $\U$}

\author{Hau-Wen Huang}
\address[H.-W. Huang]{
Department of Mathematics\\
National Central University\\
Chung-Li 32001 Taiwan}
\email{hauwenh@math.ncu.edu.tw}

\author{Chia-Yi Wen}
\address[C.-Y. Wen]{
Department of Applied Mathematics\\
National Yang Ming Chiao Tung University\\
Hsinchu 30010 Taiwan}
\email{cywen.sc08@nycu.edu.tw}

\begin{document}

\begin{abstract}
From the viewpoint of Johnson graphs as slices of a hypercube, we derive a novel algebra homomorphism $\sharp$ from the universal Racah algebra $\Re$ into $\U$.
We use the Casimir elements of $\Re$ to describe the kernel of $\sharp$.
By pulling back via $\sharp$ every $\U$-module can be viewed as an $\Re$-module.
We show that for any finite-dimensional $\U$-module $V$, the $\Re$-module $V$ is completely reducible and three generators of $\Re$
act on every irreducible $\Re$-submodule of $V$ as a Leonard triple.
In particular, Leonard triples can be constructed in terms of the second dual distance operator of the hypercube $H(D,2)$ and a decomposition of the second distance operator of $H(D,2)$ induced by Johnson graphs.
\end{abstract}

\maketitle



\allowdisplaybreaks

\section{Introduction}\label{sec:intro}

Throughout this paper we adopt the following conventions:
A vacuous summation is defined as zero. A vacuous product is defined as one. Let $\N$ denote the set of all nonnegative integers.
Let $\Z$ denote the ring of integers.
Let $\C$ denote the complex number field and $\i$ represents the imaginary unit $\sqrt{-1}$. An {\it algebra} is meant to be a unital associative algebra. A {\it subalgebra} has the same unit as the parent algebra. A subalgebra is said to be {\it proper} if the subalgebra is not equal to the parent algebra and is not generated by the unit. An {\it algebra homomorphism} is meant to be a unital algebra homomorphism. In an algebra the notation
$[x,y]$ stands for the commutator
$
xy-yx
$.
Given any linear map $f$ the notation $\ker f$ denotes the kernel of $f$ and ${\rm Im} f$ denotes the image of $f$.
Given any set $X$ the notation $\C^X$ stands for the free vector space over $\C$ that has basis $X$.
Given any vector space $V$ a direct sum of $n$ copies of $V$ is denoted by $n\cdot V$ for any integer $n\geq 1$.

\begin{defn}\label{defn:U}
The {\it universal enveloping algebra $\U$ of $\mathfrak{sl}_2$} is an algebra over $\C$ generated by $E,F,H$ subject to the relations
\begin{align}
& [H,E]=2E, \label{eq:U_relation_E}
\\
& [H,F]=-2F, \label{eq:U_relation_F}
\\
& [E,F]=H. \label{eq:U_relation_H}
\end{align}
\end{defn}
The {\it Casimir element} $\Lambda$ of $\U$ is a remarkable central element of $\U$ defined by
\begin{equation}\label{U_casimir}
  \Lambda=EF+FE+\frac{H^2}{2}.
\end{equation}

Recall from \cite[Section 1.2]{Koe2010} that the Racah polynomials are a family of discrete orthogonal polynomials and can be expressed in terms of the hypergeometric series ${}_4F_3$. The Racah polynomials are the most general hypergeometric orthogonal polynomials in the Askey scheme.
The Racah algebra is an algebra defined by generators and relations that demonstrates the bispectral properties of Racah polynomials. Originally the Racah algebra appeared in the study of the quantum mechanical coupling of three angular momenta \cite{Levy1965}. Two decades later the Racah algebra rematerialized in the study of the symmetries of the $6j$-symbol \cite{zhedanov1988}.
Over the years the Racah algebra has developed from the initial application to become an independent subject \cite{SH:2019-1,SH:2017-1,PSarah2024}. Additionally it has been associated with the superintegrable systems
\cite{Racah&Howedual:2018,gvz2014,R&LD2014,Kalnins1,Kalnins2,R&models2019},
the additive double affine Hecke algebra of type $(C_1^\vee,C_1)$ \cite{R&BI2015,Huang:R<BI,Huang:R<BImodules},
the algebra $\U$
\cite{gvz2013,SH:2019-2}
and Leonard pairs \cite{Vidunas:2007,Askeyscheme}.
In this paper we work on a central extension of the Racah algebra:

\begin{defn}
[Equations (2.8) and (2.9), \cite{Levy1965}]
\label{defn:Re}
The {\it universal Racah algebra} $\Re$ is an algebra over $\C$ generated by $A,B,C,\Delta$ and the relations assert that
\begin{align}\label{R0}
[A,B]=[B,C]=[C,A]=2\Delta
\end{align}
and each of the following elements is central in $\Re$:
\begin{align}
&[A,\Delta]+AC-BA,
\label{R1}
\\
&[B,\Delta]+BA-CB,
\label{R2}
\\
&[C,\Delta]+CB-AC.
\label{R3}
\end{align}
\end{defn}

Using (\ref{R0}) yields that $\Re$ is generated by $A,B,C$ and
\begin{align}\label{R4}
 A+B+C
\end{align}
is central in $\Re$. For convenience let $\alpha,\beta,\gamma,\delta$ denote the central elements (\ref{R1})--(\ref{R4}) of $\Re$ respectively.
From the viewpoint of Johnson graphs as slices of a hypercube, we come up with a new connection between $\U$ and $\Re$:

\begin{thm}\label{thm:RetoU}
There exists a unique algebra homomorphism $\sharp:\Re\rightarrow \U$ that sends
\begin{eqnarray}
A &\mapsto & \frac{(E+F-2)(E+F+2)}{16}, \label{A}
\\
B & \mapsto & \frac{(H-2)(H+2)}{16}, \label{B}
\\
C & \mapsto & \frac{(\i E-\i F-2)(\i E-\i F+2)}{16}, \label{C}
\\
\Delta &\mapsto &  \frac{(H+2)F^2-(H-2)E^2}{64}. \label{D}
\end{eqnarray}
Moreover $\sharp$ sends
\begin{eqnarray}\label{sharp:alpha-delta}
\alpha &\mapsto &0,
\qquad
\beta \;\;\mapsto\; \; 0,
\qquad
\gamma \;\;\mapsto\; \; 0,
\qquad
\delta \;\;\mapsto\; \; \frac{\Lambda-6}{8}.
\end{eqnarray}
\end{thm}

Let $\mathfrak C$ denote the subalgebra of $\Re$ generated by $\alpha,\beta,\gamma,\delta$. The {\it Casimir class} of $\Re$ is an individual coset of $\mathfrak C$ in $\Re$ \cite[Section 6]{SH:2017-1}. Each element of the Casimir class is central in $\Re$ and is called a {\it Casimir element} of $\Re$.
The algebra $\Re$ contains the following symmetric Casimir elements:
\begin{align}
\Omega_A&=\Delta^2+\frac{BAC+CAB}{2}+A^2+B\gamma-C\beta-A\delta, \label{Omega_A}
\\
\Omega_B&=\Delta^2+\frac{CBA+ABC}{2}+B^2+C\alpha-A\gamma-B\delta, \label{Omega_B}
\\
\Omega_C&=\Delta^2+\frac{ACB+BCA}{2}+C^2+A\beta-B\alpha-C\delta. \label{Omega_C}
\end{align}

\begin{thm}
\label{thm:image_of_omega}
The images of $\Omega_A,\Omega_B,\Omega_C$ under $\sharp$ are equal to
\begin{gather}
\label{Omega_sharp}
-\frac{3}{1024}(\Lambda-4)(\Lambda-12).
\end{gather}
\end{thm}

\noindent Furthermore we determine the generators of $\ker \sharp$:

\begin{thm}
\label{thm:kernel_of_sharp}
Suppose that $\Omega\in \{\Omega_A,\Omega_B,\Omega_C\}$.
Then
\begin{gather*}
\alpha,
\qquad
\beta,
\qquad
256 \Omega+3(4\delta-3)(4\delta+1)
\end{gather*}
form a minimal generating set for the two-sided ideal $\ker \sharp$ of $\Re$.
\end{thm}

By pulling back via $\sharp$, each $\U$-module can be regarded as an $\Re$-module.
The next part of this paper is to study finite-dimensional $\U$-modules as $\Re$-modules.
Our investigation reveals that all these $\Re$-modules are related to Leonard triples. The notion of Leonard triples \cite{cur2007-2} was initially motivated by the notable connection among spin models, distance-regular graphs and Leonard pairs  \cite{Jaeger1992,Nomura1999,cur2007-1}. Thenceforth, Leonard triples have gradually become an independent topic \cite{Huang:2012,LTRacah,LTBannai/Ito:odd,LTBannai/Ito:even} and been found connections to other areas such as the Terwilliger algebras of some $Q$-polynomial distance-regular graphs \cite{odd:2024,Lthypercube,LTRacah&DRG}, the Racah--Wigner coefficients of $U_q(\mathfrak{sl}_2)$ \cite{Huang:RW}, the anitcommutator spin algebra \cite{antispin} and the double affine Hecke algebra of type $(C_1^\vee,C_1)$ \cite{Huang:DAHA&LT}.
The definition of Leonard triples is presented as below.
A square matrix is called {\it tridiagonal} whenever each nonzero entry is on the diagonal, the subdiagonal, or the superdiagonal. A tridiagonal matrix is called {\it irreducible} whenever each entry on the subdiagonal is nonzero, and each entry on the superdiagonal is nonzero.

\begin{defn}
[Definition 1.2, \cite{cur2007-2}]
\label{defn:Leonard_triple}
Let $V$ denote a vector space with finite positive dimension. An ordered triple of linear operators $L:V\rightarrow V,L^*:V\rightarrow V,L^\e:V\rightarrow V$ is called a {\it Leonard triple} if the following conditions hold:
\begin{enumerate}
  \item There exists a basis for $V$ with respect to which the matrix representing $L$ is diagonal and the matrices representing $L^*$ and $L^\e$ are irreducible tridiagonal.

  \item There exists a basis for $V$ with respect to which the matrix representing $L^*$ is diagonal and the matrices representing $L^\e$ and $L$ are irreducible tridiagonal.
  \item There exists a basis for $V$ with respect to which the matrix representing $L^\e$ is diagonal and the matrices representing $L$ and $L^*$ are irreducible tridiagonal.
\end{enumerate}
\end{defn}


Let $\mathcal A$ denote an algebra.
Recall that a nonzero $\mathcal A$-module is said to be {\it irreducible} whenever it has exactly two $\mathcal A$-submodules, namely $\{0\}$ and itself. An $\mathcal A$-module is {\it completely reducible} whenever it is isomorphic to a direct sum of irreducible $\mathcal A$-modules.

\begin{thm}
\label{thm:R-module decomposition}
Suppose that $V$ is a finite-dimensional $\U$-module. Then the following statements hold:
\begin{enumerate}
\item The $\Re$-module $V$ is completely reducible.

\item $A, B, C$ act on every irreducible $\Re$-submodule of $V$ as a Leonard triple.
\end{enumerate}
\end{thm}

Fix an integer $D\geq 2$. Let $X$ denote the power set of a $D$-element set.
The {\it $D$-dimensional hypercube} $H(D,2)$ is a finite simple connected graph with vertex set $X$ and two vertices $x,y\in X$ are adjacent if and only if $x$ covers $y$ or $y$ covers $x$ in the poset $(X,\subseteq)$.
The second distance relation $R_2$ of $H(D,2)$ is related to the Johnson graphs in the following way: Let $k$ denote an integer with $0\leq k\leq D$.
The {\it Johnson graph} $J(D,k)$ is a finite simple connected graph with vertex set consisting of all $x\in X$ with $|x|=k$ and two vertices $x,y$ of $J(D,k)$ are adjacent if and only if $(x,y)\in R_2$. This induces a decomposition of the second distance operator of $H(D,2)$ into the sum of the two linear maps $\A_2^J:\C^{X}\to \C^{X}$ and $\A_2^{\bar J}:\C^{X}\to \C^{X}$  given by
\begin{align}
\A_2^J x=\sum_{\substack{|y|=|x|\\(x,y)\in R_2}} y
\qquad
\hbox{for all $x\in X$},
\label{eq:A2+}
\\
\A_2^{\bar J} x=\sum_{\substack{|y|\not=|x|\\(x,y)\in R_2}} y
\qquad
\hbox{for all $x\in X$}.
\label{eq:A2-}
\end{align}
Note that $\A_2^J$ is a direct sum of the adjacency operators of $J(D,k)$ for all $k=0,1,\ldots, D$.
The vertex set of $J(D,k)$ is the set of all vertices of $H(D,2)$ at distance $k$ from the empty set $\emptyset$.
By \cite{BannaiIto1984,TerAlgebraI,TerAlgebraII,TerAlgebraIII} the second dual distance operator $\A^*_2$ of $H(D,2)$ with respect to $\emptyset$ is a linear map $\C^{X}\to \C^{X}$ given by
\begin{equation}\label{eq:dual_A_hypercube}
\A^*_2 x=\frac{(D-2|x|)^2-D}{2} x
\qquad
\hbox{for all $x\in X$}.
\end{equation}
Let $\R$ denote the subalgebra of ${\rm End}(\C^{X})$ generated by $\A_2^J,\A_2^{\bar J},\A_2^{*}$. The matrices representing $\A_2^J,\A_2^{\bar J},\A_2^{*}$ with respect to the basis $X$ for $\C^{X}$ are invariant under conjugate-transpose. Therefore the finite-dimensional algebra $\R$ is semisimple and every irreducible $\R$-module is isomorphic to an irreducible $\R$-submodule of the standard $\R$-module $\C^{X}$.

We obtain from Theorems \ref{thm:RetoU} and \ref{thm:R-module decomposition} a connection among the algebra $\R$, the universal Racah algebra $\Re$ and Leonard triples:

\begin{thm}\label{thm:LP_hypercube+johnson}
{ }
\begin{enumerate}
\item There exists a unique algebra homomorphism $\Re\to \R$ that sends
\begin{eqnarray}
A &\mapsto & \frac{D}{16}-\frac{1}{4}+\frac{\A_2^J+\A_2^{\bar J}}{8},
\label{A:R->R}
\\
B &\mapsto & \frac{D}{16}-\frac{1}{4}+\frac{\A_2^*}{8},
\label{B:R->R}
\\
C &\mapsto & \frac{D}{16}-\frac{1}{4}+\frac{\A_2^J-\A_2^{\bar J}}{8}.
\label{C:R->R}
\end{eqnarray}
Moreover the above algebra homomorphism is surjective.

\item
$\A_2^J+\A_2^{\bar J}$,
$\A_2^{*}$,
$\A_2^J-\A_2^{\bar J}$ act on each irreducible $\R$-module as a Leonard triple.
\end{enumerate}
\end{thm}

The paper is organized as follows:
In \S\ref{sec:D_3-action on R and U(sl2)} we give a $D_3$-action on $\U$ and recall a $D_3$-action on $\Re$ from \cite{SH:2017-1}.
In \S\ref{sec:RetoU} we prove Theorem \ref{thm:RetoU} and display that $\sharp$ makes a connection between the $D_3$-actions on $\U$ and $\Re$.
In \S\ref{sec:image_of_omega} we prove Theorem \ref{thm:image_of_omega} by applying the $\Z$-graded algebra structure of $\U$ and the even subalgebra $\U_e$ of $\U$.
In \S\ref{sec:image of omega and kernel of sharp} we prove Theorem \ref{thm:kernel_of_sharp} by employing a basis for $\Re$ given in \cite{SH:2017-1} and a basis for $\U_e$ given in \cite{halved:2023}.
In \S\ref{sec:Remodule} we recall finite-dimensional irreducible $\Re$-modules from \cite{SH:2019-1} and
establish the necessary and sufficient conditions for $A,B,C$ acting on finite-dimensional irreducible $\Re$-modules as Leonard triples.
In \S\ref{sec:Remodule_Ln} we show Theorem~\ref{thm:R-module decomposition} with the help of the results stated in \S\ref{sec:Remodule}.
In \S\ref{sec:Proof_for_thm_LP_hypercube+johnson} we show Theorem \ref{thm:LP_hypercube+johnson} by integrating Theorems \ref{thm:RetoU} and \ref{thm:R-module decomposition} along with a $\U$-module structure on $\C^{X}$ given by \cite{hypercube2002}.
In \S\ref{sec:R_structure} we finish this paper by studying the algebraic properties of $\R$.

\section{The $D_3$-actions on $\U$ and $\Re$}
\label{sec:D_3-action on R and U(sl2)}

The dihedral group $D_n$ for any integer $n\geq 3$ is the group of symmetries of a regular $n$-sided polygon.
Recall that $D_3$ has a presentation with the generators $\sigma,\tau$ and the relations
\begin{gather}\label{rel:D3}
\sigma^2=1,
\qquad
\tau^3=1,
\qquad
(\sigma\tau)^2=1.
\end{gather}

\begin{prop}
\label{prop:D3onU}
There exists a unique $D_3$-action on $\U$ such that each of the following holds:
\begin{enumerate}
\item $\sigma$ acts on $\U$ as an algebra automorphism in the following way:

\begin{table}[H]
\centering
\extrarowheight=3pt
\begin{tabular}{c|ccc}
$u$  &$E$ &$F$ &$H$

\\

\midrule[1pt]

$\sigma(u)$
&$\i F$
&$-\i E$
&$-H$

\end{tabular}
\end{table}

\item $\tau$ acts on $\U$ as an algebra automorphism in the following way:

\begin{table}[H]
\centering
\extrarowheight=3pt
\begin{tabular}{c|ccc}
$u$  &$E$ &$F$ &$H$

\\

\midrule[1pt]

$\tau(u)$
&$\displaystyle\frac{H-\i E-\i F}{2}$
&$\displaystyle\frac{H+\i E+\i F}{2}$
&$\i E-\i F$
\end{tabular}
\end{table}
\end{enumerate}
\end{prop}
\begin{proof}
Using Definition~\ref{defn:U} it is routine to verify that there exists a unique algebra automorphism $\mathcal{S}$ of $\U$ that sends
\begin{eqnarray*}
(E,F,H)
&\mapsto &
\left(
\i F,-\i E,-H
\right)
\end{eqnarray*}
and there exists a unique algebra automorphism $\mathcal{T}$ of $\U$ that sends
\begin{eqnarray*}
(E,F,H)
&\mapsto &
\left(
\frac{H-\i E-\i F}{2},\frac{H+\i E+\i F}{2},\i E-\i F
\right).
\end{eqnarray*}
The composition $\mathcal{S}\circ \mathcal{T}$ sends
\begin{eqnarray*}
(E,F,H) &\mapsto & \left(-\frac{H+E-F}{2}, -\frac{H-E+F}{2}, -E-F\right).
\end{eqnarray*}
It is straightforward to check that
$\mathcal S^2=1$, $\mathcal T^3=1$ and $(\mathcal S\circ \mathcal T)^2=1$.
By (\ref{rel:D3}) there exists a unique $D_3$-action on $\U$ such that $\sigma$ and $\tau$ act as $\mathcal{S}$ and $\mathcal{T}$  respectively. The proposition follows.
\end{proof}

\begin{prop}
The $D_3$-action on $\U$ is faithful.
\end{prop}
\begin{proof}
The group $D_3$ consists of the following elements:
\begin{gather}
\label{D3_U}
1,
\qquad
\tau,
\qquad
\tau^2,
\qquad
\sigma,
\qquad
\sigma \tau,
\qquad
\sigma \tau^2.
\end{gather}
By Proposition \ref{prop:D3onU} the algebra automorphisms (\ref{D3_U}) of $\U$ send $H$ to
\begin{gather}
\label{D3_H}
H,
\qquad
\i E-\i F,
\qquad
E+F,
\qquad
-H,
\qquad
-E-F,
\qquad
\i F-\i E,
\end{gather}
respectively.

Let $x$ denote an indeterminate over $\C$. By Definition \ref{defn:U} there exists a unique $\U$-module structure on a polynomial ring $\C[x]$ given by
\begin{align*}
E x^i &=
\left\{
\begin{array}{ll}
(1-i) x^{i-1}
\qquad &\hbox{if $i\geq 1$},
\\
0
\qquad &\hbox{if $i=0$}
\end{array}
\right.
\qquad
\hbox{for all $i\in \N$},
\\
F x^i&=(i+1) x^{i+1}
\qquad
\hbox{for all $i\in \N$},
\\
H x^i&=-2i x^i
\qquad
\hbox{for all $i\in \N$}.
\end{align*}
Since the elements (\ref{D3_H}) are mutually distinct operators on the $\U$-module $\C[x]$, the algebra automorphisms (\ref{D3_U}) of $\U$ are mutually distinct. The proposition follows.
\end{proof}

Recall the Casimir element $\Lambda$ of $\U$ from (\ref{U_casimir}).

\begin{lem}
\label{lem:D3&U_casimr}
$g(\Lambda)=\Lambda$ for all $g\in D_3$.
\end{lem}
\begin{proof}
By Proposition~\ref{prop:D3onU}(i) the algebra automorphism $\sigma$ of $\U$ maps
\begin{eqnarray*}
  EF &\mapsto & FE \;\;\mapsto \;\; EF
  ,\qquad
  H^2 \;\;\mapsto \;\; H^2.
\end{eqnarray*}
Combined with (\ref{U_casimir}) this implies that $\sigma(\Lambda)=\Lambda$.
By Proposition~\ref{prop:D3onU}(ii) the algebra automorphism $\tau$ of $\U$ maps
\begin{eqnarray*}
  EF &\mapsto &\frac{H^2
  +\i[H,E+F]
  +(E+F)^2}{4}, \\
  FE &\mapsto & \frac{H^2
  -\i[H,E+F]
  +(E+F)^2}{4}, \\
  H^2 &\mapsto & -(E-F)^2.
\end{eqnarray*}
Combined with (\ref{U_casimir}) this implies that $\tau(\Lambda)=\Lambda$.
The lemma follows since the group $D_3$ is generated by $\sigma$ and $\tau$.
\end{proof}

The following $D_3$-action on $\Re$ is obtained from the $D_6$-action on $\Re$ given in \cite[Proposition 4.1]{SH:2017-1}.

\begin{prop}
[Proposition 4.1, \cite{SH:2017-1}]
\label{prop:D3onRe}
There exists a unique $D_3$-action on $\Re$ such that each of the following holds:
\begin{enumerate}
\item $\sigma$ acts on $\Re$ as an algebra automorphism in the following way:

\begin{table}[H]
\centering
\extrarowheight=3pt
\begin{tabular}{c|cccc|cccc}
$u$  &$A$ &$B$ &$C$ &$\Delta$
&$\alpha$ &$\beta$ &$\gamma$ &$\delta$
\\

\midrule[1pt]

$\sigma(u)$ &$C$ &$B$ &$A$ &$-\Delta$
&$-\gamma$ &$-\beta$ &$-\alpha$ &$\delta$
\end{tabular}
\end{table}

\item $\tau$ acts on $\Re$ as an algebra automorphism  in the following way:

\begin{table}[H]
\centering
\extrarowheight=3pt
\begin{tabular}{c|cccc|cccc}
$u$  &$A$ &$B$ &$C$ &$\Delta$
&$\alpha$ &$\beta$ &$\gamma$ &$\delta$
\\

\midrule[1pt]

$\tau(u)$ &$B$ &$C$ &$A$ &$\Delta$
&$\beta$ &$\gamma$ &$\alpha$ &$\delta$
\end{tabular}
\end{table}
\end{enumerate}
\end{prop}

\begin{prop}
[Proposition 4.3, \cite{SH:2017-1}]
The $D_3$-action on $\Re$ is faithful.
\end{prop}

Recall the Casimir elements $\Omega_A,\Omega_B,\Omega_C$ of $\Re$ from (\ref{Omega_A})--(\ref{Omega_C}).

\begin{lem}
[Lemma 6.3, \cite{SH:2017-1}]
\label{lem:D3&Re_casimir}
The set $\{\Omega_A,\Omega_B,\Omega_C\}$ is invariant under the $D_3$-action on $\Re$. Moreover the $D_3$-action on $\{\Omega_A,\Omega_B,\Omega_C\}$ is as follows:

\begin{table}[H]
\centering
\extrarowheight=3pt
\begin{tabular}{c|ccc}
$u$  &$\Omega_A$ &$\Omega_B$ &$\Omega_C$
\\

\midrule[1pt]

$\sigma(u)$ &$\Omega_C$ &$\Omega_B$ &$\Omega_A$
\\
$\tau(u)$ &$\Omega_B$ &$\Omega_C$ &$\Omega_A$
\end{tabular}
\end{table}
\end{lem}

In the subsequent section we will see how the $D_3$-action on $\U$ is related to the $D_3$-action on $\Re$.

\section{Proof for Theorem \ref{thm:RetoU}}
\label{sec:RetoU}

Observe that
\begin{equation}\label{eq:U_identity}
 [x,yz]=y[x,z]+[x,y]z
\end{equation}
for any elements $x,y,z$ in an algebra.

\begin{lem}
\label{lem:H^2E,H^2F}
The following equations hold in $\U$:
  \begin{enumerate}
    \item $[H^2,E]=4(H-1)E$.
    \item $[H^2,F]=-4(H+1)F$.
    \item $[H^2,E^2]=8(H-2)E^2$.
    \item $[H^2,F^2]=-8(H+2)F^2$.
  \end{enumerate}
\end{lem}
\begin{proof}
(i): The equation (i) follows by applying (\ref{eq:U_identity}) with $(x,y,z)=(E,H,H)$ and simplifying the resulting equation by using (\ref{eq:U_relation_E}).

(iii):  The equation (iii) follows by applying (\ref{eq:U_identity}) with $(x,y,z)=(H^2,E,E)$ and simplifying the resulting equation by using (\ref{eq:U_relation_E}) and Lemma~\ref{lem:H^2E,H^2F}(i).

 (ii), (iv):  By Proposition~\ref{prop:D3onU}(i) the equations (ii) and (iv) follow by applying $\sigma$ to  Lemma~\ref{lem:H^2E,H^2F}(i), (iii).
\end{proof}

\noindent {\it Proof of Theorem \ref{thm:RetoU}.}
Let $A^{\sharp},B^{\sharp},C^{\sharp},\Delta^{\sharp}$ denote the right-hand sides of (\ref{A})--(\ref{D}) respectively.
To see the existence of $\sharp$, by Definition \ref{defn:Re} it suffices to prove that the following equations hold in $\U$:
\begin{gather}
  [A^\sharp,B^\sharp]=[B^\sharp,C^\sharp]=[C^\sharp,A^\sharp]=2\Delta^\sharp,
  \label{R0:sharp}\\
  [A^\sharp,\Delta^\sharp]+A^\sharp C^\sharp-B^\sharp A^\sharp=0,
  \label{R1:sharp}\\
  [B^\sharp,\Delta^\sharp]+B^\sharp A^\sharp-C^\sharp B^\sharp=0,
  \label{R2:sharp}\\
  [C^\sharp,\Delta^\sharp]+C^\sharp B^\sharp-A^\sharp C^\sharp=0.
  \label{R3:sharp}
\end{gather}

Expanding $A^{\sharp},B^{\sharp},C^{\sharp}$ directly and by (\ref{U_casimir}) it follows that
\begin{gather}
\label{R4:sharp}
 A^\sharp+B^\sharp+C^\sharp=\frac{\Lambda-6}{8}.
 \end{gather}
Since $\Lambda$ is central in $\U$ the equation (\ref{R4:sharp}) implies the first two equalities in (\ref{R0:sharp}).
By Lemma \ref{lem:H^2E,H^2F}(i), (ii) the commutator $[H^2,E+F]$ is equal to $4$ times
$(H-1)E-(H+1)F$.
Combined with (\ref{eq:U_relation_E})--(\ref{eq:U_relation_H}) this yields that
\begin{align*}
\frac{[H^2,E+F](E+F)}{4}&=H(H+1)+(H-1)E^2-(H+1)F^2-2EF,
\\
\frac{(E+F)[H^2,E+F]}{4}&=-H(H+1)+(H-3)E^2-(H+3)F^2+2EF.
\end{align*}
The relation $[A^\sharp, B^\sharp]=2 \Delta^\sharp$ follows by applying (\ref{eq:U_identity}) with $(x,y,z)=(H^2,E+F-2,E+F+2)$ and using the above two equations to simplify the resulting equation. Therefore (\ref{R0:sharp}) holds.

A direct expansion yields that
\begin{align}
B^{\sharp}A^{\sharp}&=\frac{(H^2-4)(E^2+F^2+EF+FE-4)}{256},
\label{BA:sharp}\\
B^{\sharp}C^{\sharp}&=-\frac{(H^2-4)(E^2+F^2-EF-FE+4)}{256}.
\notag
\end{align}
By (\ref{R0:sharp}) the element $C^{\sharp}B^{\sharp}=B^{\sharp}C^{\sharp}-2 \Delta^\sharp$. Hence
\begin{align}
\label{CB:sharp}
C^{\sharp}B^{\sharp}=\frac{(H-2)E^2-(H+2)F^2}{32}-\frac{(H^2-4)(E^2+F^2-EF-FE+4)}{256}.
\end{align}
Applying (\ref{eq:U_identity}) with $(x,y,z)=(H^2,H-2,E^2)$, $(H^2,H+2,F^2)$ we evaluate the resulting equations by using  Lemma \ref{lem:H^2E,H^2F}(iii), (iv). It follows that
\begin{align*}
[H^2,(H-2)E^2]&=8(H-2)^2 E^2,
\\
[H^2,(H+2)F^2]&=-8(H+2)^2 F^2.
\end{align*}
Hence
\begin{gather}
\label{[BD]:sharp}
[B^{\sharp},\Delta^{\sharp}]=-\frac{(H-2)^2 E^2+(H+2)^2 F^2}{128}.
\end{gather}
Now (\ref{R2:sharp}) follows by using (\ref{BA:sharp})--(\ref{[BD]:sharp}) to evaluate the left-hand side of (\ref{R2:sharp}). By Proposition \ref{prop:D3onU}(ii) the algebra automorphism $\tau$ of $\U$ maps
\begin{eqnarray*}
E+F  &\mapsto & H \;\; \mapsto \; \; \i E-\i F\;\; \mapsto \; \; E+F.
\end{eqnarray*}
Hence $\tau$ maps
\begin{eqnarray*}
A^\sharp  &\mapsto & B^\sharp \;\; \mapsto \; \; C^\sharp \;\; \mapsto \; \; A^\sharp.
\end{eqnarray*}
Combined with (\ref{R0:sharp}) the automorphism $\tau$ fixes $\Delta^\sharp$.
The equations (\ref{R1:sharp}) and (\ref{R3:sharp}) follow by applying $\tau^{-1}$ and $\tau$ to (\ref{R2:sharp}) respectively.

By (\ref{R0:sharp})--(\ref{R3:sharp}) the existence of $\sharp$ follows. By Definition \ref{defn:Re} the algebra $\Re$ is generated by $A,B,C,\Delta$. The uniqueness of $\sharp$ follows.
By (\ref{R1:sharp})--(\ref{R3:sharp}) the map $\sharp$ sends each of $\alpha,\beta,\gamma$ to zero. By (\ref{R4:sharp}) the image of $\delta$ under $\sharp$ is equal to $\frac{\Lambda-6}{8}$. The result follows.
\hfill $\square$

\medskip

Recall the $D_3$-actions on $\U$ and $\Re$ from Propositions \ref{prop:D3onU} and \ref{prop:D3onRe} respectively.

\begin{thm}\label{thm:D3_commutes}
For any $g\in D_3$ the following diagram commutes:
\begin{table}[H]
\centering
\begin{tikzpicture}
\matrix(m)[matrix of math nodes,
row sep=4em, column sep=6em,
text height=1.5ex, text depth=0.25ex]
{
\Re
&\U\\
\Re
&\U\\
};
\path[->,font=\scriptsize,>=angle 90]
(m-1-1) edge node[above] {$\sharp$} (m-1-2)
(m-2-1) edge node[below] {$\sharp$} (m-2-2)
(m-1-1) edge node[left] {$g$} (m-2-1)
(m-1-2) edge node[right] {$g$} (m-2-2);
\end{tikzpicture}
\end{table}
\end{thm}
\begin{proof}
In the proof of Theorem \ref{thm:RetoU} we have seen that the algebra automorphism $\tau$ of $\U$ cyclically permutes the right-hand sides of (\ref{A})--(\ref{C}). Combined with Proposition \ref{prop:D3onRe}(ii) the algebra homomorphisms $\tau\circ\sharp$ and $\sharp\circ \tau$ agree at $A,B,C$. By Proposition \ref{prop:D3onU}(i) the algebra automorphism $\sigma$ of $\U$ maps
\begin{eqnarray*}
E+F  &\mapsto & \i F-\i E \;\; \mapsto \; \; E+F,
\qquad
H  \;\; \mapsto \; \; -H.
\end{eqnarray*}
Hence the algebra automorphism $\sigma$ of $\U$ fixes (\ref{B}) and swaps (\ref{A}) and (\ref{C}).
Combined with Proposition \ref{prop:D3onRe}(i) the algebra homomorphisms $\sigma\circ\sharp$ and $\sharp\circ \sigma$ agree at $A,B,C$.

By (\ref{R0}) the algebra $\Re$ is generated by $A,B,C$. Therefore the diagram commutes for $g\in\{\sigma,\tau\}$. The result follows since the group $D_3$ is generated by $\sigma$ and $\tau$.
\end{proof}

\section{Proof for Theorem \ref{thm:image_of_omega} }
\label{sec:image_of_omega}


Recall that an algebra $\mathcal A$ is called a {\it $\Z$-graded algebra} whenever there are subspaces $\{\mathcal A_n\}_{n\in \Z}$ of $\mathcal A$ satisfying the following properties:
\begin{description}

\item[(G1)] $\mathcal A=\bigoplus\limits_{n\in\Z}\mathcal A_n$.

\item[(G2)] $\mathcal A_m\cdot \mathcal A_n\subseteq \mathcal A_{m+n}$ for all $m,n\in\Z$.

\end{description}
For any $n\in \Z$ the space $\mathcal A_n$ is called the {\it $n^{ th}$ homogeneous subspace} of $\mathcal A$. By {\bf (G1)}, for any $u\in \mathcal A$ there are unique $u_n\in \mathcal A_n$ for all $n\in \Z$ such that
$$
u=\sum_{n\in \Z} u_n.
$$
For each $n\in \Z$ the element $u_n$ is called the {\it$n^{ th}$ homogeneous component} of $u$. For any $n\in \Z$ let $U_n$ denote the subspace of $\U$ spanned by
$$
  E^iF^jH^k\qquad\hbox{for all $i,j,k\in\N$ with $i-j=n$}.
$$
It is well known that $\U$ is a $\Z$-graded algebra by setting $U_n$ as the $n^{\rm th}$ homogeneous subspace of $\U$ for each $n\in \Z$. Note that $E\in U_1$, $F\in U_{-1}$, $H\in U_0$ and $\Lambda\in U_0$.

\begin{defn}
[Definition 1.4, \cite{halved:2023}]
\label{defn:Ue}
Define
$$
\U_e=\bigoplus\limits_{n\in\Z}
U_{2n}.
$$
Since $1\in U_0$ and $\U_e$ is closed under the multiplication of $\U$, it follows that $\U_e$ is a subalgebra of $\U$. We call $\U_e$ the {\it even subalgebra} of $\U$.
\end{defn}

\begin{prop}
\label{prop:sharp&Ue}
${\rm Im}\, \sharp$ is a subalgebra of $\U_e$.
\end{prop}
\begin{proof}
Using {\bf (G2)} yields that the right-hand sides of (\ref{A})--(\ref{D}) are in
$$U_{-2}\oplus U_0\oplus U_2,
\quad
U_0,
\quad
U_{-2}\oplus U_0\oplus U_2,
\quad
U_{-2}\oplus U_2,
$$
respectively. Combined with Definitions \ref{defn:Re} and \ref{defn:Ue} the result follows.
\end{proof}

\begin{lem}
[Lemma 3.2,\cite{halved:2023}]
\label{lem:Ue_basis}
For all $n\in \N$ the following statements hold:
\begin{enumerate}
\item The elements
$$
E^{2n}\Lambda^i H^k
\qquad\hbox{for all $i,k\in\N$}
$$
are a basis for $U_{2n}$.

\item The elements
$$
F^{2n}\Lambda^iH^k
\qquad\hbox{for all $i,k\in\N$}
$$
are a basis for $U_{-2n}$.
\end{enumerate}
\end{lem}


\begin{lem}
\label{lem:Ue_relation}
The following equations hold in $\U_e$:
\begin{enumerate}
\item $[H,E^2]=4 E^2$.

\item $[H,F^2]=-4 F^2$.

\item $16E^2F^2=(H^2-2H-2\Lambda)(H^2-6H-2\Lambda+8)$.

\item $16F^2E^2=(H^2+2H-2\Lambda)(H^2+6H-2\Lambda+8)$.

\end{enumerate}
\end{lem}
\begin{proof}
Immediate from \cite[Lemmas 2.2 and 2.4]{halved:2023}.
\end{proof}


For notational convenience let $u^\sharp$ denote the image of $u$ under $\sharp$ for any $u\in \Re$.

\begin{lem}
\label{lem:sharp&ABCD}
The following equations hold in $\U_e$:
\begin{enumerate}
\item $A^\sharp=\ds\frac{\Lambda+E^2+F^2}{16}-\frac{H^2}{32}-\frac{1}{4}$.

\item $B^\sharp=\ds\frac{H^2}{16}-\frac{1}{4}$.

\item $C^\sharp=\ds\frac{\Lambda-E^2-F^2}{16}-\frac{H^2}{32}-\frac{1}{4}$.

\item $\Delta^\sharp=\ds\frac{F^2(H-2)-E^2(H+2)}{64}$.
\end{enumerate}
\end{lem}
\begin{proof}
The equation (ii) is immediate from (\ref{B}).
The equations (i) and (iii) follow by using (\ref{U_casimir}) to rewrite the right-hand sides of (\ref{A}) and (\ref{C}) respectively. The equation (iv) follows by applying Lemma \ref{lem:Ue_relation}(i), (ii) to rewrite the right-hand side of (\ref{D}).
\end{proof}

\begin{lem}
\label{lem:Re_homogenous}
All nonzero homogeneous components of $A^\sharp,B^\sharp, C^\sharp, \Delta^\sharp, \delta^\sharp$ are as follows:

\begin{table}[H]
\centering
\extrarowheight=3pt
\begin{tabular}{c||c|c|c|c|c}
$u$  &$A^\sharp$ &$B^\sharp$ &$C^\sharp$ &$\Delta^\sharp$ &$\delta^\sharp$
\\

\hline \hline
$u_{-2}$
&$\frac{F^2}{16}$
&$0$
&$-\frac{F^2}{16}$
&$\frac{F^2H}{64}-\frac{F^2}{32}$
&$0$
\\
\hline
$u_0$
&$\frac{\Lambda}{16}-\frac{H^2}{32}-\frac{1}{4}$
&$\frac{H^2}{16}-\frac{1}{4}$
&$\frac{\Lambda}{16}-\frac{H^2}{32}-\frac{1}{4}$
&$0$
&$\frac{\Lambda-6}{8}$
\\
\hline
$u_2$
&$\frac{E^2}{16}$
&$0$
&$-\frac{E^2}{16}$
&$-\frac{E^2H}{64}-\frac{E^2}{32}$
&$0$
\end{tabular}
\end{table}
\end{lem}
\begin{proof}
By Lemma \ref{lem:Ue_basis} the nonzero homogeneous components of $A^\sharp,B^\sharp, C^\sharp, \Delta^\sharp$ are
immediately obtained from Lemma \ref{lem:sharp&ABCD}.
By (\ref{sharp:alpha-delta}) the element $\delta^\sharp=\frac{\Lambda-6}{8}\in U_0$.
\end{proof}

\noindent{\it Proof of Theorem~\ref{thm:image_of_omega}.}
Recall the Casimir element $\Omega_B$ of $\Re$ from (\ref{Omega_B}).
By Lemma \ref{lem:Re_homogenous} and since $\alpha^\sharp=\gamma^\sharp=0$ by (\ref{sharp:alpha-delta}), the element $\Omega_B^\sharp\in U_{-4}\oplus U_{-2}\oplus U_0\oplus U_2\oplus U_4$ and
\begin{align*}
{(\Omega_B^\sharp)}_{-4}
&=(\Delta^\sharp_{-2})^2
+
\frac{C^\sharp_{-2} B^\sharp_0 A^\sharp_{-2}
+
A^\sharp_{-2} B^\sharp_0 C^\sharp_{-2}}{2},
\\
{(\Omega_B^\sharp)}_{-2}&=
\frac{C^\sharp_{-2} B^\sharp_0 A^\sharp_0
+
C^\sharp_0 B^\sharp_0 A^\sharp_{-2}
+
A^\sharp_{-2} B^\sharp_0 C^\sharp_0
+
A^\sharp_0 B^\sharp_0 C^\sharp_{-2}}{2},
\\
{(\Omega_B^\sharp)}_0&=
\Delta^\sharp_2 \Delta^\sharp_{-2}
+
\Delta^\sharp_{-2}\Delta^\sharp_2
+
(B^\sharp_0)^2
-
B^\sharp_0\delta^\sharp_0
+
A^\sharp_0 B^\sharp_0 C^\sharp_0
\\
&\qquad +\,
\frac{C^\sharp_2 B^\sharp_0 A^\sharp_{-2}
+
C^\sharp_{-2} B^\sharp_0 A^\sharp_2
+
A^\sharp_2 B^\sharp_0 C^\sharp_{-2}
+
A^\sharp_{-2} B^\sharp_0 C^\sharp_2}{2}.
\end{align*}

Concerning the equation for ${(\Omega_B^\sharp)}_{-4}$, it follows from Lemmas \ref{lem:Ue_relation}(ii) and  \ref{lem:Re_homogenous} that only the terms in the first row of the table below may have nonzero coefficients in $(\Delta^\sharp_{-2})^2$, $C^\sharp_{-2} B^\sharp_0 A^\sharp_{-2}$, $A^\sharp_{-2} B^\sharp_0 C^\sharp_{-2}$ with respect to the basis for $U_{-4}$ given in Lemma \ref{lem:Ue_basis}(ii).
A routine calculation shows that those coefficients are as follows:
\begin{table}[H]
\centering
\extrarowheight=3pt
\begin{tabular}{c||c|c|c}
coefficients  &$F^4$ &$F^4 H$ &$F^4 H^2$
\\
\hline \hline
$(\Delta^\sharp_{-2})^2$
&$\frac{3}{1024}$
&$-\frac{1}{512}$
&$\frac{1}{4096}$
\\
\hline
$C^\sharp_{-2} B^\sharp_0 A^\sharp_{-2}$
&$-\frac{3}{1024}$
&$\frac{1}{512}$
&$-\frac{1}{4096}$
\\
\hline
$A^\sharp_{-2} B^\sharp_0 C^\sharp_{-2}$
&$-\frac{3}{1024}$
&$\frac{1}{512}$
&$-\frac{1}{4096}$
\end{tabular}
\end{table}
\noindent From the above table we can conclude that ${(\Omega_B^\sharp)}_{-4}=0$.

Concerning the equation for ${(\Omega_B^\sharp)}_{-2}$, it follows from Lemmas \ref{lem:Ue_relation}(ii) and  \ref{lem:Re_homogenous} that only the terms in the first row of the table below
may have nonzero coefficients in $C^\sharp_{-2} B^\sharp_0 A^\sharp_0$, $C^\sharp_0 B^\sharp_0 A^\sharp_{-2}$,
$A^\sharp_{-2} B^\sharp_0 C^\sharp_0$, $A^\sharp_0 B^\sharp_0 C^\sharp_{-2}$ with respect to the basis for $U_{-2}$ given in Lemma \ref{lem:Ue_basis}(ii).
A routine calculation shows that those coefficients are as follows:
\begin{table}[H]
\centering
\extrarowheight=3pt
\begin{tabular}{c||c|c|c|c|c|c|c|c}
coefficients
&$F^2 \Lambda$
&$F^2 \Lambda H$
&$F^2 \Lambda H^2$
&$F^2$
&$F^2H$
&$F^2H^2$
&$F^2H^3$
&$F^2H^4$
\\

\hline \hline

$C^\sharp_{-2} B^\sharp_0 A^\sharp_0$
&$\frac{1}{1024}$
&$0$
&$-\frac{1}{4096}$
&$-\frac{1}{256}$
&$0$
&$\frac{1}{2048}$
&$0$
&$\frac{1}{8192}$
\\
\hline

$C^\sharp_0 B^\sharp_0 A^\sharp_{-2}$
&$\frac{3}{1024}$
&$-\frac{1}{512}$
&$\frac{1}{4096}$
&$-\frac{9}{256}$
&$\frac{9}{256}$
&$-\frac{25}{2048}$
&$\frac{1}{512}$
&$-\frac{1}{8192}$
\\
\hline

$A^\sharp_{-2} B^\sharp_0 C^\sharp_0$
&$-\frac{1}{1024}$
&$0$
&$\frac{1}{4096}$
&$\frac{1}{256}$
&$0$
&$-\frac{1}{2048}$
&$0$
&$-\frac{1}{8192}$
\\
\hline

$A^\sharp_0 B^\sharp_0 C^\sharp_{-2}$
&$-\frac{3}{1024}$
&$\frac{1}{512}$
&$-\frac{1}{4096}$
&$\frac{9}{256}$
&$-\frac{9}{256}$
&$\frac{25}{2048}$
&$-\frac{1}{512}$
&$\frac{1}{8192}$
\end{tabular}
\end{table}
\noindent From the above table we can conclude that ${(\Omega_B^\sharp)}_{-2}=0$.

Concerning the equation for ${(\Omega_B^\sharp)}_0$, it follows from Lemmas \ref{lem:Ue_relation} and  \ref{lem:Re_homogenous} that only the terms in the first rows of the tables below may have nonzero coefficients in $\Delta^\sharp_2 \Delta^\sharp_{-2}$,
$\Delta^\sharp_{-2}\Delta^\sharp_2$,
$(B^\sharp_0)^2$,
$B^\sharp_0\delta^\sharp_0$,
$A^\sharp_0 B^\sharp_0 C^\sharp_0$,
$C^\sharp_2 B^\sharp_0 A^\sharp_{-2}$,
$C^\sharp_{-2} B^\sharp_0 A^\sharp_2$,
$A^\sharp_2 B^\sharp_0 C^\sharp_{-2}$,
$A^\sharp_{-2} B^\sharp_0 C^{\sharp}_2$ with respect to the basis for $U_0$ given in Lemma \ref{lem:Ue_basis}.
A routine but tedious calculation shows that those coefficients are as follows:
\begin{table}[H]
\centering
\extrarowheight=3pt
\begin{tabular}{c||c|c|c|c|c|c|c|c}
coefficients
&$\Lambda^2$ &$\Lambda^2 H$ &$\Lambda^2 H^2$
&$\Lambda$ &$\Lambda H$ &$\Lambda H^2$
&$\Lambda H^3$ &$\Lambda H^4$
\\

\hline \hline

$\Delta^\sharp_2 \Delta^\sharp_{-2}$
&$-\frac{1}{4096}$
&$\frac{1}{4096}$
&$-\frac{1}{16384}$
&$\frac{1}{1024}$
&$-\frac{1}{512}$
&$\frac{3}{2048}$
&$-\frac{1}{2048}$
&$\frac{1}{16384}$
\\
\hline
$\Delta^\sharp_{-2}\Delta^\sharp_2$
&$-\frac{1}{4096}$
&$-\frac{1}{4096}$
&$-\frac{1}{16384}$
&$\frac{1}{1024}$
&$\frac{1}{512}$
&$\frac{3}{2048}$
&$\frac{1}{2048}$
&$\frac{1}{16384}$
\\
\hline
$(B^\sharp_0)^2$
&$0$
&$0$
&$0$
&$0$
&$0$
&$0$
&$0$
&$0$
\\
\hline
$B^\sharp_0\delta^\sharp_0$
&$0$
&$0$
&$0$
&$-\frac{1}{32}$
&$0$
&$\frac{1}{128}$
&$0$
&$0$
\\
\hline
$A^\sharp_0 B^\sharp_0 C^\sharp_0$
&$-\frac{1}{1024}$
&$0$
&$\frac{1}{4096}$
&$\frac{1}{128}$
&$0$
&$-\frac{1}{1024}$
&$0$
&$-\frac{1}{4096}$
\\
\hline
$C^\sharp_2 B^\sharp_0 A^\sharp_{-2}$
&$-\frac{3}{4096}$
&$\frac{1}{2048}$
&$-\frac{1}{16384}$
&$\frac{3}{1024}$
&$-\frac{5}{1024}$
&$\frac{3}{1024}$
&$-\frac{3}{4096}$
&$\frac{1}{16384}$
\\
\hline
$C^\sharp_{-2} B^\sharp_0 A^\sharp_2$
&$-\frac{3}{4096}$
&$-\frac{1}{2048}$
&$-\frac{1}{16384}$
&$\frac{3}{1024}$
&$\frac{5}{1024}$
&$\frac{3}{1024}$
&$\frac{3}{4096}$
&$\frac{1}{16384}$
\\
\hline
$A^\sharp_2 B^\sharp_0 C^\sharp_{-2}$
&$-\frac{3}{4096}$
&$\frac{1}{2048}$
&$-\frac{1}{16384}$
&$\frac{3}{1024}$
&$-\frac{5}{1024}$
&$\frac{3}{1024}$
&$-\frac{3}{4096}$
&$\frac{1}{16384}$
\\
\hline
$A^\sharp_{-2} B^\sharp_0 C^{\sharp}_2$
&$-\frac{3}{4096}$
&$-\frac{1}{2048}$
&$-\frac{1}{16384}$
&$\frac{3}{1024}$
&$\frac{5}{1024}$
&$\frac{3}{1024}$
&$\frac{3}{4096}$
&$\frac{1}{16384}$
\end{tabular}
\end{table}

\begin{table}[H]
\centering
\extrarowheight=3pt
\begin{tabular}{c||c|c|c|c|c|c|c}
coefficients
&$1$ &$H$ &$H^2$ &$H^3$ &$H^4$ &$H^5$ &$H^6$
\\

\hline \hline

$\Delta^\sharp_2 \Delta^\sharp_{-2}$
&$0$
&$\frac{1}{1024}$
&$-\frac{9}{4096}$
&$\frac{1}{512}$
&$-\frac{7}{8192}$
&$\frac{3}{16384}$
&$-\frac{1}{65536}$
\\
\hline
$\Delta^\sharp_{-2}\Delta^\sharp_2$
&$0$
&$-\frac{1}{1024}$
&$-\frac{9}{4096}$
&$-\frac{1}{512}$
&$-\frac{7}{8192}$
&$-\frac{3}{16384}$
&$-\frac{1}{65536}$
\\
\hline
$(B^\sharp_0)^2$
&$\frac{1}{16}$
&$0$
&$-\frac{1}{32}$
&$0$
&$\frac{1}{256}$
&$0$
&$0$
\\
\hline
$B^\sharp_0\delta^\sharp_0$
&$\frac{3}{16}$
&$0$
&$-\frac{3}{64}$
&$0$
&$0$
&$0$
&$0$
\\
\hline
$A^\sharp_0 B^\sharp_0 C^\sharp_0$
&$-\frac{1}{64}$
&$0$
&$0$
&$0$
&$\frac{3}{4096}$
&$0$
&$\frac{1}{16384}$
\\
\hline
$C^\sharp_2 B^\sharp_0 A^\sharp_{-2}$
&$0$
&$\frac{3}{1024}$
&$-\frac{23}{4096}$
&$\frac{17}{4096}$
&$-\frac{3}{2048}$
&$\frac{1}{4096}$
&$-\frac{1}{65536}$
\\
\hline
$C^\sharp_{-2} B^\sharp_0 A^\sharp_2$
&$0$
&$-\frac{3}{1024}$
&$-\frac{23}{4096}$
&$-\frac{17}{4096}$
&$-\frac{3}{2048}$
&$-\frac{1}{4096}$
&$-\frac{1}{65536}$
\\
\hline
$A^\sharp_2 B^\sharp_0 C^\sharp_{-2}$
&$0$
&$\frac{3}{1024}$
&$-\frac{23}{4096}$
&$\frac{17}{4096}$
&$-\frac{3}{2048}$
&$\frac{1}{4096}$
&$-\frac{1}{65536}$
\\
\hline
$A^\sharp_{-2} B^\sharp_0 C^{\sharp}_2$
&$0$
&$-\frac{3}{1024}$
&$-\frac{23}{4096}$
&$-\frac{17}{4096}$
&$-\frac{3}{2048}$
&$-\frac{1}{4096}$
&$-\frac{1}{65536}$
\end{tabular}
\end{table}
\noindent From the above two tables we can conclude that ${(\Omega_B^\sharp)}_0$ is equal to (\ref{Omega_sharp}).

Applying Proposition~\ref{prop:D3onU}(i) yields that
\begin{gather}
\label{tsUn}
\sigma(U_n)=U_{-n}
\qquad
\hbox{for all $n\in \Z$}.
\end{gather}
By Lemma \ref{lem:D3&Re_casimir} the algebra automorphism $\sigma$ of $\Re$ fixes $\Omega_B$. Combined with Theorem \ref{thm:D3_commutes} this implies that the algebra automorphism $\sigma$ of $\U$ fixes $\Omega_B^
\sharp$. Applying (\ref{tsUn}) yields that
\begin{align*}
{(\Omega_B^\sharp)}_{4}=\sigma({(\Omega_B^\sharp)}_{-4})=0,
\\
{(\Omega_B^\sharp)}_{2}=\sigma({(\Omega_B^\sharp)}_{-2})=0.
\end{align*}
Therefore $\Omega_B^\sharp$ is equal to (\ref{Omega_sharp}).
By Lemma~\ref{lem:D3&U_casimr} the algebra automorphism $\tau$ of $\U$ fixes $\Lambda$.
By Lemma~\ref{lem:D3&Re_casimir} the algebra automorphism $\tau$ of $\Re$ maps
$$
\Omega_B \; \; \mapsto \; \; \Omega_C \; \; \mapsto \; \; \Omega_A.
$$
Combined with Theorem~\ref{thm:D3_commutes} this yields that $\Omega_A^\sharp$ and $\Omega_C^\sharp$ are equal to  (\ref{Omega_sharp}).
\hfill $\square$

\section{Proof for Theorem \ref{thm:kernel_of_sharp} }
\label{sec:image of omega and kernel of sharp}

Throughout this section we adopt the following notations: Let $x,y$ denote two commuting indeterminates over $\C$ and set
$$
f(x,y)=256 x+3(4y-3)(4y+1).
$$
Let $\Omega$ stand for any of $\Omega_A,\Omega_B,\Omega_C$.

By (\ref{sharp:alpha-delta}) the element $\delta^\sharp=\frac{\Lambda-6}{8}$. By Theorem \ref{thm:image_of_omega} the element $\Omega^\sharp=-\frac{3}{1024}(\Lambda-4)(\Lambda-12)$.
It is routine to verify that
\begin{gather}
\label{f(Omega,delta)}
f(\Omega,\delta)^\sharp=0.
\end{gather}

\begin{lem}
\label{lem:Omega&delta}
For any $p(x,y)\in \C[x,y]$ the following conditions are equivalent:
\begin{enumerate}
\item $p(\Omega,\delta)^\sharp=0$.

\item $f(x,y)$ divides $p(x,y)$.
\end{enumerate}
\end{lem}
\begin{proof}
We view any element of $\C[x,y]$ as a polynomial in one variable $x$ over $\C[y]$. Let $p(x,y)\in \C[x,y]$ be given.
Since the leading coefficient of $f(x,y)$ is invertible in $\C[y]$, there are $q(x,y)\in \C[x,y]$ and $r(y)\in \C[y]$ such that
$$
p(x,y)=q(x,y)f(x,y)+r(y).
$$
Combined with (\ref{f(Omega,delta)}) the condition (i) holds if and only if
\begin{gather}
\label{r(delta)}
r(\delta)^\sharp=0.
\end{gather}
By Lemma \ref{lem:Ue_basis} the element $\Lambda$ is algebraically independent over $\C$ and so is $\delta^\sharp$. Hence (\ref{r(delta)}) is equivalent to $r(y)=0$. The lemma follows.
\end{proof}

\begin{lem}
\label{lem:kernel}
For all $n\in \N$ the following statements hold:
\begin{enumerate}
\item The elements
$$
E^{2n} H^k
\qquad
\hbox{for all $k\in \N$}
$$
are a basis for $U_{2n}$ over $\C[\Lambda]$.

\item The elements
$$
F^{2n} H^k
\qquad
\hbox{for all $k\in \N$}
$$
are a basis for $U_{-2n}$ over $\C[\Lambda]$.
\end{enumerate}
\end{lem}
\begin{proof}
Immediate from Lemma \ref{lem:Ue_basis}.
\end{proof}

\begin{lem}
\label{lem:Re_basis}
The elements
$$
A^i \Delta^j B^k \Omega^\ell \alpha^r \delta^s \beta^t
\qquad
\hbox{for all $i,k,\ell,r,s,t\in \N$ and $j\in \{0,1\}$}
$$
are a basis for $\Re$.
\end{lem}
\begin{proof}
Immediate from \cite[Lemma 7.1 and Theorem 7.5]{SH:2017-1}.
\end{proof}

\begin{prop}
\label{prop:kernel_of_sharp}
$\ker \sharp$ is the
two-sided ideal of $\Re$ generated by
$\alpha, \beta, f(\Omega,\delta)$.
\end{prop}
\begin{proof}
By (\ref{sharp:alpha-delta}) the elements $\alpha$ and $\beta$ are in $\ker \sharp$. By (\ref{f(Omega,delta)}) the element $f(\Omega,\delta)$ is in $\ker \sharp$. Hence $\ker \sharp$ contains the two-sided ideal $\I$ of $\Re$ generated by
$\alpha$, $\beta$ and $f(\Omega,\delta)$.
Let $u$ denote any element of $\ker \sharp$.
It suffices to show that $u\in \I$.

Recall that $\alpha$ and $\beta$ are central in $\Re$.
Combined with Lemma \ref{lem:Re_basis}, there are $v,w\in \Re$ such that $u-v\alpha-w\beta$ is a linear combination of
$$
A^i \Delta^j B^k \Omega^\ell \delta^s
\qquad
\hbox{for all $i,k,\ell,s\in \N$ and $j\in \{0,1\}$}.
$$
In other words, there are $n\in \N$ and $p_{ik}(x,y),q_{ik}(x,y)\in \C[x,y]$ for all integers $i,k$ with $0\leq i, k\leq n$ such that
\begin{equation}\label{eq:u}
u=
v\alpha
+w\beta
+
\sum_{i=0}^n \sum_{k=0}^n A^i B^k p_{ik}(\Omega,\delta)+\sum_{i=0}^n \sum_{k=0}^n A^i \Delta B^k q_{ik}(\Omega,\delta).
\end{equation}
Since $\alpha^\sharp=0$ and $\beta^\sharp=0$ it follows that
\begin{gather}
\label{e1:usharp}
u^\sharp
=
\sum_{i=0}^n \sum_{k=0}^n
{(A^\sharp)}^i {(B^\sharp)}^k p_{ik}(\Omega,\delta)^\sharp
+
\sum_{i=0}^n \sum_{k=0}^n
{(A^\sharp)}^i \Delta^\sharp {(B^\sharp)}^k q_{ik}(\Omega,\delta)^\sharp.
\end{gather}

We are going to show that
\begin{gather}
\label{p=0&q=0}
p_{ik}(\Omega,\delta)^\sharp=q_{ik}(\Omega,\delta)^\sharp=0
\qquad
(0\leq i, k\leq n).
\end{gather}
Suppose on the contrary that there exists an integer $i$ with $0\leq i\leq n$ such that
\begin{align*}
&\hbox{$p_{ik}(\Omega,\delta)^\sharp\not=0$
for some integer $k$ with $0\leq k\leq n$}
\\
\hbox{or }&\hbox{$q_{ik}(\Omega,\delta)^\sharp\not=0$ for some integer $k$ with $0\leq k\leq n$}.
\end{align*}
Let $\ell$ denote the largest integer $i$ satisfying the above condition.
The equation (\ref{e1:usharp}) now can be written as
\begin{gather}
\label{e2:usharp}
u^\sharp
=
\sum_{i=0}^\ell \sum_{k=0}^n
{(A^\sharp)}^i {(B^\sharp)}^k p_{ik}(\Omega,\delta)^\sharp
+
\sum_{i=0}^\ell \sum_{k=0}^n
{(A^\sharp)}^i \Delta^\sharp {(B^\sharp)}^k q_{ik}(\Omega,\delta)^\sharp.
\end{gather}
Applying Lemma \ref{lem:Re_homogenous} yields that $u^\sharp_{2\ell+2}$ is equal to
\begin{gather*}
{(A^\sharp_2)}^\ell \Delta^\sharp_2
\sum_{k=0}^n {(B^\sharp_0)}^k q_{\ell k}(\Omega,\delta)^\sharp
=
-\frac{E^{2\ell+2}(H+2)}{4^{2\ell+3}}
\sum_{k=0}^n
\frac{(H^2-4)^k}{16^{k}} q_{\ell k}(\Omega,\delta)^\sharp.
\end{gather*}
On the other hand, since $u^\sharp=0$ it follows that $u^\sharp_{2\ell+2}=0$.
Hence
\begin{gather}
\label{e1:u2l+2}
E^{2\ell+2}(H+2)
\sum_{k=0}^n
\frac{(H^2-4)^k}{16^{k}} q_{\ell k}(\Omega,\delta)^\sharp
=0.
\end{gather}
In what follows we will apply Lemma~\ref{lem:kernel}(i) without further mention.
By looking at the coefficients of $E^{2\ell+2}H^{2k+1}$ in both sides of (\ref{e1:u2l+2}) in the ordering $k=n,n-1,\ldots,0$, it can be deduced that
\begin{gather}
\label{qlk=0}
q_{\ell k}(\Omega,\delta)^\sharp=0
\qquad
(0\leq k\leq n).
\end{gather}

By (\ref{qlk=0}) the equation (\ref{e2:usharp}) can be simplified to
\begin{gather*}
u^\sharp
=
\sum_{i=0}^\ell \sum_{k=0}^n
{(A^\sharp)}^i {(B^\sharp)}^k p_{ik}(\Omega,\delta)^\sharp
+
\sum_{i=0}^{\ell-1} \sum_{k=0}^n
{(A^\sharp)}^i \Delta^\sharp {(B^\sharp)}^k q_{ik}(\Omega,\delta)^\sharp.
\end{gather*}
The second double summation is interpreted as zero if $\ell=0$.
Applying Lemma \ref{lem:Re_homogenous} yields that $u^\sharp_{2\ell}$ is equal to
\begin{gather*}
\frac{E^{2\ell}}{16^\ell}
\sum_{k=0}^n \frac{(H^2-4)^k}{16^k}
p_{\ell k}(\Omega,\delta)^\sharp
-
\frac{E^{2\ell}(H+2)}{4^{2\ell+1}}
\sum_{k=0}^n
\frac{(H^2-4)^k}{16^k}
q_{\ell-1,k}(\Omega,\delta)^\sharp.
\end{gather*}
If $\ell=0$ then $q_{\ell-1,k}(\Omega,\delta)^\sharp$ is interpreted as zero for any integer $k$ with $0\leq k\leq n$.
Since $u^\sharp_{2\ell}=0$ it follows that
\begin{gather}
\label{e1:u2l}
4 E^{2\ell}
\sum_{k=0}^n \frac{(H^2-4)^k}{16^k}
p_{\ell k}(\Omega,\delta)^\sharp
=
E^{2\ell}(H+2)
\sum_{k=0}^n
\frac{(H^2-4)^k}{16^k}
q_{\ell-1,k}(\Omega,\delta)^\sharp.
\end{gather}
By looking at the coefficients of $E^{2\ell}H^{2k+1}$ in both sides of (\ref{e1:u2l}) in the ordering $k=n,n-1,\ldots,0$, it can be deduced that
\begin{gather*}
q_{\ell-1, k}(\Omega,\delta)^\sharp=0
\qquad
(0\leq k\leq n)
\end{gather*}
even if $\ell\geq 1$. Hence the right-hand side of (\ref{e1:u2l}) is equal to zero. Now looking at the coefficients of $E^{2\ell}H^{2k}$ in both sides of (\ref{e1:u2l}) in the ordering $k=n,n-1,\ldots,0$ yields that
\begin{gather}
\label{plk=0}
p_{\ell k}(\Omega,\delta)^\sharp=0
\qquad
(0\leq k\leq n).
\end{gather}
The results (\ref{qlk=0}) and (\ref{plk=0}) contradict to the setting of $\ell$. Therefore (\ref{p=0&q=0}) is indeed true. By Lemma \ref{lem:Omega&delta} the polynomial $f(x,y)$ is a common divisor of $p_{ik}(x,y)$ and $q_{ik}(x,y)$ for all integers $i,k$ with $0\leq i,k\leq n$.  Combined with (\ref{eq:u}) this implies that $u\in \I$. The result follows.
\end{proof}

\begin{lem}
\label{lem:minimal}
\begin{enumerate}
\item $\alpha$
is not in the two-sided ideal of $\Re$ generated by $f(\Omega,\delta)$ and $\beta$.
\item $\beta$
is not in the two-sided ideal of $\Re$ generated by $f(\Omega,\delta)$ and $\alpha$.
\item $f(\Omega,\delta)$ is not in the two-sided ideal of $\Re$ generated by $\alpha$ and $\beta$.
\end{enumerate}
\end{lem}
\begin{proof}
(i): Suppose on the contrary that there are $v,w\in \Re$ such that
\begin{gather}
\label{minimal_2}
\alpha=
v f(\Omega,\delta)+ w \beta.
\end{gather}
By Lemma \ref{lem:Re_basis} we may assume that $v$ and $w$ are linear combinations of
\begin{align*}
&A^i \Delta^j B^k \Omega^\ell \alpha^r \delta^s
\qquad
\hbox{for all $i,k,\ell,r,s\in \N$ and $j\in \{0,1\}$},
\\
&A^i \Delta^j B^k \Omega^\ell \alpha^r \delta^s \beta^t
\qquad
\hbox{for all $i,k,\ell,r,s,t\in \N$ and $j\in \{0,1\}$},
\end{align*}
respectively.
Clearly the coefficients of
\begin{align*}
A^i \Delta^j B^k \Omega^{\ell+1} \alpha^r \delta^s,
\qquad
A^i \Delta^j B^k \Omega^\ell \alpha^r \delta^s \beta^{t+1}
\end{align*}
in the left-hand side of (\ref{minimal_2}) with respect to the basis for $\Re$ mentioned in Lemma \ref{lem:Re_basis} are equal to zero for all $i,k,\ell,r,s,t\in \N$ and $j\in \{0,1\}$. Hence $v=0$ and $w=0$.
Then the right-hand side of (\ref{minimal_2}) is equal to zero, a contradiction.

(ii), (iii): Similar to the proof for Lemma \ref{lem:minimal}(i).
\end{proof}

\noindent{\it Proof of Theorem~\ref{thm:kernel_of_sharp}.}
Immediate from Proposition \ref{prop:kernel_of_sharp} and Lemma \ref{lem:minimal}.
\hfill $\square$


\section{Finite-dimensional irreducible $\Re$-modules and Leonard triples}
\label{sec:Remodule}

\begin{prop}
[Proposition 2.4, \cite{SH:2019-1}]
\label{prop:Rd(abc)}
For any $a,b,c\in \C$ and $d\in \N$ there exists a $(d+1)$-dimensional $\Re$-module $R_d(a,b,c)$ satisfying the following conditions:
\begin{enumerate}
  \item There exists a basis for $R_d(a,b,c)$ with respect to which the matrices representing $A$ and $B$ are
\begin{equation*}
  \begin{pmatrix}
    \theta_0 &  &  &  &  \\
    1 & \theta_1 &  &  &  \\
     & 1 & \theta_2 &  &  \\
     &  & \ddots & \ddots &  \\
     &  &  & 1 & \theta_d
  \end{pmatrix},\qquad
\begin{pmatrix}
    \theta_0^* & \varphi_1 &  &  &  \\
     & \theta_1^* & \varphi_2 &  &  \\
     &  & \theta_2^* & \ddots &  \\
     &  &  & \ddots & \varphi_d \\
     &  &  &  & \theta_d^*
  \end{pmatrix},
\end{equation*}
respectively, where
\begin{align*}
  \theta_i & =(a+\tfrac{d}{2}-i)(a+\tfrac{d}{2}-i+1)\qquad
  (0\leq i\leq d),
 \\
  \theta_i^* & =(b+\tfrac{d}{2}-i)(b+\tfrac{d}{2}-i+1)\qquad
  (0\leq i\leq d),
  \\
  \varphi_i &=i(i-d-1)(a+b+c+\tfrac{d}{2}-i+2)(a+b-c+\tfrac{d}{2}-i+1)\qquad
  (1\leq i\leq d).
  \notag
\end{align*}
  \item The elements $\alpha,\beta,\gamma,\delta$ act on $R_d(a,b,c)$ as scalar multiplication by
\begin{align*}
&(c-b)(c+b+1)(a-\tfrac{d}{2})(a+\tfrac{d}{2}+1),
\\
&(a-c)(a+c+1)(b-\tfrac{d}{2})(b+\tfrac{d}{2}+1),
\\
&(b-a)(b+a+1)(c-\tfrac{d}{2})(c+\tfrac{d}{2}+1),
\\
&\tfrac{d}{2}(\tfrac{d}{2}+1)+a(a+1)+b(b+1)+c(c+1),
\end{align*}
respectively.
\end{enumerate}
\end{prop}

\begin{lem}
\label{lem:trace_Rd(abc)}
For any $a,b,c\in \C$ and $d\in \N$ the traces of $A,B,C$ on the $\Re$-module $R_d(a,b,c)$ are equal to $d+1$ times
\begin{align*}
a(a+1)+\frac{d(d+2)}{12},
\qquad
b(b+1)+\frac{d(d+2)}{12},
\qquad
c(c+1)+\frac{d(d+2)}{12},
\end{align*}
respectively.
\end{lem}
\begin{proof}
It is straightforward to evaluate the traces of $A$ and $B$ on the $\Re$-module $R_d(a,b,c)$ by applying Proposition \ref{prop:Rd(abc)}(i). Recall that $\delta$ is equal to (\ref{R4}). Using this relation and Proposition \ref{prop:Rd(abc)}(ii) it is easy to evaluate the trace of $C$ on $R_d(a,b,c)$.
\end{proof}

\begin{lem}
[Theorem 4.5, \cite{SH:2019-1}]
\label{lem:irr_Rd(a,b,c)}
For any $a,b,c\in \C$ and $d\in \N$ the following conditions are equivalent:
\begin{enumerate}
\item The $\Re$-module $R_d(a,b,c)$ is irreducible.

\item $
a+b+c+1,-a+b+c,a-b+c,a+b-c
\not\in\{\tfrac{d}{2}-i\,|\,i=1,2,\ldots,d\}.
$
\end{enumerate}
\end{lem}

According to \cite[Theorem 6.3]{SH:2019-1}, for any finite-dimensional irreducible $\Re$-module $V$ there are $a,b,c\in \C$ and $d\in \N$ such that the $\Re$-module $R_d(a,b,c)$ is isomorphic to $V$. The number $d$ is equal to the dimension of $V$ minus one. Furthermore the parameters $a,b,c$ can be determined in the following way:

\begin{lem}
[Corollary 6.5, \cite{SH:2019-1}]
\label{lem:Rd(abc)_isomorphism class}
Suppose that $V$ is a $(d+1)$-dimensional irreducible $\Re$-module. For any $a,b,c\in \C$ the following conditions are equivalent:
\begin{enumerate}
\item The $\Re$-module $R_d(a,b,c)$ is isomorphic to $V$.

\item The traces of $A,B,C$ on the $\Re$-module $V$ are equal to $d+1$ times
\begin{align*}
  a(a+1)+\frac{d(d+2)}{12},
  \qquad
  b(b+1)+\frac{d(d+2)}{12},
  \qquad
  c(c+1)+\frac{d(d+2)}{12},
\end{align*}
respectively.
\end{enumerate}
\end{lem}

Let $V$ denote any $\Re$-module. Given any algebra automorphism $g$ of $\Re$ the notation
$V^g$
stands for the $\Re$-module obtained by twisting the $\Re$-module $V$ via $g$.

\begin{lem}
\label{lem:Rd(abc)_tau}
For any $a,b,c\in \C$ and $d\in \N$ the following conditions are equivalent:
\begin{enumerate}
\item The $\Re$-module $R_d(a,b,c)$ is irreducible.

\item The $\Re$-module $R_d(c,a,b)$ is irreducible.

\item The $\Re$-module $R_d(b,c,a)$ is irreducible.
\end{enumerate}
Suppose that {\rm (i)--(iii)} hold. Then the $\Re$-modules
$R_d(a,b,c)$, $R_d(c,a,b)^\tau$, $R_d(b,c,a)^{\tau^2}$
are isomorphic.
\end{lem}
\begin{proof}
The equivalence of (i)--(iii) is immediate from Lemma \ref{lem:irr_Rd(a,b,c)}.
By Proposition \ref{prop:D3onRe}(ii) and Lemma \ref{lem:trace_Rd(abc)} the elements $A,B,C$ have the same traces on the $\Re$-modules
$R_d(a,b,c)$, $R_d(c,a,b)^\tau$, $R_d(b,c,a)^{\tau^2}$.
Thus, by Lemma \ref{lem:Rd(abc)_isomorphism class} these $\Re$-modules are isomorphic provided that (i)--(iii) hold.
\end{proof}


\begin{lem}
\label{lem:Rd(abc)_ABCminpoly}
Suppose that $a,b,c\in \C$ and $d\in \N$ such that the $\Re$-module $R_d(a,b,c)$ is irreducible.
Then the minimal polynomials of $A,B,C$ on the $\Re$-module $R_d(a,b,c)$ are
$$
\prod_{i=0}^d(x-\theta_i),
\qquad
\prod_{i=0}^d(x-\theta_i^*),
\qquad
\prod_{i=0}^d(x-\theta_i^\e),
$$
respectively, where $x$ denotes an indeterminate over $\C$ and
\begin{align*}
\theta_i &=\textstyle (a+\frac{d}{2}-i)(a+\frac{d}{2}-i+1)
\qquad
(0\leq i\leq d),
\\
\theta_i^* &=\textstyle (b+\frac{d}{2}-i)(b+\frac{d}{2}-i+1)
\qquad
(0\leq i\leq d),
\\
\theta_i^\e &=\textstyle (c+\frac{d}{2}-i)(c+\frac{d}{2}-i+1)
\qquad
(0\leq i\leq d).
\end{align*}
\end{lem}
\begin{proof}
By Proposition \ref{prop:Rd(abc)}(i) the minimal polynomial of $A$ on $R_d(a,b,c)$ is
$\prod_{i=0}^d (x-\theta_i)$.
By Proposition \ref{prop:D3onRe}(ii) the action of $C$ on the $\Re$-module $R_d(c,a,b)^\tau$ is identical to the action of $A$ on the $\Re$-module $R_d(c,a,b)$.
By Lemma \ref{lem:Rd(abc)_tau} the $\Re$-module $R_d(c,a,b)^\tau$ is isomorphic to $R_d(a,b,c)$.
Hence the minimal polynomial of $C$ on the $\Re$-module $R_d(a,b,c)$ is equal to the minimal polynomial of $A$ on $R_d(c,a,b)$, namely $\prod_{i=0}^d (x-\theta_i^\e)$.
By Proposition \ref{prop:D3onRe}(ii) the action of $B$ on the $\Re$-module $R_d(b,c,a)^{\tau^2}$ is identical to the action of $A$ on the $\Re$-module $R_d(b,c,a)$. By Lemma \ref{lem:Rd(abc)_tau} the $\Re$-module $R_d(b,c,a)^{\tau^2}$ is isomorphic to $R_d(a,b,c)$.
Hence the minimal polynomial of $B$ on the $\Re$-module $R_d(a,b,c)$ is equal to the minimal polynomial of $A$ on $R_d(b,c,a)$, namely $\prod_{i=0}^d (x-\theta_i^*)$.
\end{proof}

\begin{lem}
\label{lem:Rd(abc)_ABCdiag}
Suppose that $a,b,c\in \C$ and $d\in \N$ such that the $\Re$-module $R_d(a,b,c)$ is irreducible. Then the following conditions are equivalent:
\begin{enumerate}
\item $A$ {\rm(}resp. $B${\rm)} {\rm(}resp. $C${\rm)} is diagonalizable on $R_d(a,b,c)$.

\item $a$ {\rm(}resp. $b${\rm)} {\rm(}resp. $c${\rm)} is not in $\{\frac{i-d-1}{2}\,|\,i=1,2,\ldots,2d-1\}$.
\end{enumerate}
\end{lem}
\begin{proof}
It is routine to verify the lemma by using Lemma \ref{lem:Rd(abc)_ABCminpoly}.
\end{proof}

\begin{lem}\label{lem:relations_in_R}
The following relations hold in $\Re$:
\begin{enumerate}
\item $A^2B-2ABA+BA^2-2AB-2BA=2A^2-2A\delta+2\alpha$.

\item $B^2C-2BCB+CB^2-2BC-2CB=2B^2-2B\delta+2\beta$.

\item $C^2A-2CAC+AC^2-2CA-2AC=2C^2-2C\delta+2\gamma$.

\item $A^2C-2ACA+CA^2-2AC-2CA=2A^2-2A\delta-2\alpha$.

\item $B^2A-2BAB+AB^2-2BA-2AB=2B^2-2B\delta-2\beta$.

\item $C^2B-2CBC+BC^2-2CB-2BC=2C^2-2C\delta-2\gamma$.
\end{enumerate}
\end{lem}
\begin{proof}
(i): By (\ref{R0}) the relation $\Delta=\frac{[A,B]}{2}$ holds. Recall that $\alpha$ is equal to (\ref{R1}) and $\delta$ is equal to (\ref{R4}). The relation (i) follows by applying the above first and third relations to express the right-hand side of the above second relation in terms of $A,B,\delta$.

(ii), (iii): By Proposition \ref{prop:D3onRe}(ii) the relations (ii) and (iii) follow by applying $\tau$ and $\tau^2$ to Lemma \ref{lem:relations_in_R}(i) respectively.

(iv)--(vi): By Proposition \ref{prop:D3onRe}(i) the relations (iv)--(vi) follow by applying $\sigma$ to Lemma \ref{lem:relations_in_R}(iii), (ii), (i) respectively.
\end{proof}

Recall Leonard triples from Definition \ref{defn:Leonard_triple}.
We are now able to characterize some equivalent conditions for $A,B,C$ acting as a Leonard triple on a finite-dimensional irreducible $\Re$-module.

\begin{thm}
\label{thm:LT}
Suppose that $a,b,c\in \C$ and $d\in \N$ such that the $\Re$-module $R_d(a,b,c)$ is irreducible. Then the following conditions are equivalent:
\begin{enumerate}
\item $A,B,C$ act on $R_d(a,b,c)$ as a Leonard triple.

\item $A,B,C$ are diagonalizable on $R_d(a,b,c)$.

\item $a,b,c\not\in \{\frac{i-d-1}{2}\,|\,i=1,2,\ldots,2d-1\}$.
\end{enumerate}
\end{thm}
\begin{proof}
(i) $\Rightarrow$ (ii): Immediate from Definition~\ref{defn:Leonard_triple}.

(ii) $\Leftrightarrow$ (iii): Immediate from Lemma \ref{lem:Rd(abc)_ABCdiag}.

(ii), (iii) $\Rightarrow$ (i): Suppose that (ii) and (iii) hold. Let
$$
\theta_i=
\textstyle
(a+\frac{d}{2}-i)(a+\frac{d}{2}-i+1)
\qquad
(-1\leq i\leq d+1).
$$
Recall the minimal polynomial of $A$ on the $\Re$-module $R_d(a,b,c)$ from Lemma \ref{lem:Rd(abc)_ABCminpoly}.
Let $u_i$ $(0\leq i\leq d)$ denote a $\theta_i$-eigenvector of $A$ on the $\Re$-module $R_d(a,b,c)$. Note that the vectors $\{u_i\}_{i=0}^d$ are a basis for $R_d(a,b,c)$. Hence the matrix representing $A$ with respect to the basis $\{u_i\}_{i=0}^d$ for the $\Re$-module $R_d(a,b,c)$ is diagonal. Applying either side of Lemma \ref{lem:relations_in_R}(i) to $u_i$ ($0\leq i\leq d$) yields that
\begin{gather*}
(A-\theta_{i-1})(A-\theta_{i+1}) B u_i
\end{gather*}
is a scalar multiple of $u_i$. Applying either side of Lemma \ref{lem:relations_in_R}(iv) to $u_i$ ($0\leq i\leq d$) yields that
\begin{gather*}
(A-\theta_{i-1})(A-\theta_{i+1}) C u_i
\end{gather*}
is a scalar multiple of $u_i$.
By (ii) or (iii) the scalars $\{\theta_i\}_{i=0}^d$ are mutually distinct. Hence $Bu_i$ and $Cu_i$ ($1\leq i\leq d-1$) are linear combinations of $u_{i-1},u_i,u_{i+1}$. By (iii) the scalar $\theta_{-1}\not\in\{\theta_2,\theta_3,\ldots,\theta_d\}$. Hence $Bu_0$ and $Cu_0$ are linear combinations of $u_0$ and $u_1$. By (iii) the scalar $\theta_{d+1}\not\in\{\theta_0,\theta_1,\ldots,\theta_{d-2}\}$. Hence $Bu_d$ and $Cu_d$ are linear combinations of $u_{d-1}$ and $u_d$. By the above comments the matrices representing $B$ and $C$ with respect to the basis $\{u_i\}_{i=0}^d$ for the $\Re$-module $R_d(a,b,c)$ are tridiagonal. Applying the irreducibility of the $\Re$-module $R_d(a,b,c)$, a routine argument shows that the two tridiagonal matrices are irreducible.

In light of Lemma \ref{lem:relations_in_R}(ii), (v) (resp. Lemma \ref{lem:relations_in_R}(iii), (vi)) a similar argument shows that there exists a basis for the $\Re$-module $R_d(a,b,c)$ with respect to which the matrix representing $B$ (resp. $C$) is diagonal and the matrices representing $C$ and $A$ (resp. $A$ and $B$) are irreducible tridiagonal.
Therefore (i) follows.
\end{proof}


\section{Proof for Theorem \ref{thm:R-module decomposition}}
\label{sec:Remodule_Ln}

Let $n\in \N$ be given. By Definition \ref{defn:U} it is routine to verify that there exists an $(n+1)$-dimensional $\U$-module $L_n$ that has a basis $\{v_i\}_{i=0}^n$ such that
\begin{align*}
E v_i&=(n-i+1) v_{i-1}
\qquad
(1\leq i\leq n),
\qquad
E v_0=0,
\\
F v_i&=(i+1) v_{i+1}
\qquad
(0\leq i\leq n-1),
\qquad
F v_n=0,
\\
H v_i&=(n-2i) v_i
\qquad
(0\leq i\leq n).
\end{align*}
The $\U$-module $L_n$ is irreducible and every $(n+1)$-dimensional irreducible $\U$-module is isomorphic to $L_n$ \cite[Proposition V.4.2 and Theorem V.4.4]{kassel}. Furthermore every finite-dimensional $\U$-module is completely reducible \cite[Theorem V.4.6]{kassel}.

Since $\U_e$ is a subalgebra of $\U$ every $\U$-module can be regarded as a $\U_e$-module. Let $V$ denote a given $\U$-module.
For any $\theta\in \C$ let $V(\theta)$ denote the subspace of $V$ consisting of all $v\in V$ with $Hv=\theta v$.
Recall from \cite[Proposition 5.1]{halved:2023} that
$$
\bigoplus_{i\in \Z} V(\theta+4i)
$$
is a $\U_e$-submodule of $V$ for any $\theta\in \C$. Define
\begin{align*}
L_n^{(0)} &= \bigoplus_{i=0}^{\floor{\frac{n}{2}}} L_n(n-4i)
\qquad
\hbox{for all $n\in \N$},
\\
L_n^{(1)} &= \bigoplus_{i=0}^{\floor{\frac{n-1}{2}}} L_n(n-4i-2)
\qquad
\hbox{for all integers $n\geq 1$}.
\end{align*}
Observe that
\begin{equation}\label{eq:Ln_decomposition_as_U_e}
L_n=
\left\{
\begin{array}{ll}
L_n^{(0)}
\qquad &\hbox{if $n=0$},
\\
L_n^{(0)}\oplus L_n^{(1)}
 \qquad &\hbox{if $n\geq 1$}.
\end{array}
\right.
\end{equation}
Therefore $L_n^{(0)}$ is a $\U_e$-submodule of $L_n$ for each $n\in \N$ and $L_n^{(1)}$ is a $\U_e$-submodule of $L_n$ for each integer $n\geq 1$. By \cite[Lemma 5.5]{halved:2023} the $\U_e$-module $L_n^{(0)}$ is irreducible for each $n\in \N$. By \cite[Lemma 5.8]{halved:2023} the $\U_e$-module $L_n^{(1)}$ is irreducible for each integer $n\geq 1$.
Furthermore the $\U_{e}$-modules $L_n^{(0)}$ for all $n\in\N$ and the $\U_{e}$-modules $L_n^{(1)}$ for all integers $n\geq 1$ are mutually non-isomorphic \cite[Theorem 5.10]{halved:2023}.

\begin{lem}
[Lemma 5.4, \cite{halved:2023}]
\label{lem:Ln0_Ue}
For any $n\in \N$ the $\U_e$-module $L_n^{(0)}$ satisfies the following conditions:
\begin{enumerate}
\item There exists a basis $\{u_i^{(0)}\}_{i=0}^{\floor{\frac{n}{2}}}$ for the $\U_e$-module $L_n^{(0)}$ such that
        \begin{align*}
             E^2 u_{i}^{(0)} &=(n-2i+1)(n-2i+2) u_{i-1}^{(0)}
            \quad
            (1\leq i\leq \textstyle{\floor{\frac{n}{2}}}),
            \qquad
            E^2 u_0^{(0)}=0,
            \\
             F^2 u_{i}^{(0)} &=(2i+1)(2i+2) u_{i+1}^{(0)}
            \quad
            (0\leq i\leq \textstyle\floor{\frac{n}{2}}-1),
            \qquad
            F^2 u_{\floor{\frac{n}{2}}}^{(0)}=0,
            \\
             H u_{i}^{(0)} &=(n-4i) u_{i}^{(0)}
            \quad
            (0\leq i\leq \textstyle\floor{\frac{n}{2}}).
        \end{align*}
\item The element $\Lambda$ acts on $L_n^{(0)}$ as scalar multiplication by $\frac{n(n+2)}{2}$.
\end{enumerate}
\end{lem}

The notation $\{u_i^{(0)}\}_{i=0}^{\floor{\frac{n}{2}}}$ will always stand for the basis for the $\U_e$-module $L_n^{(0)}$ ($n\in \N$) given in Lemma \ref{lem:Ln0_Ue}.
By Proposition \ref{prop:sharp&Ue} the $\U_e$-module $L_n^{(0)}$ is an $\Re$-module by pulling back via $\sharp$.

\begin{prop}\label{prop:A,B,C on Ln0}
Suppose that $n\in \N$. Then the actions of $A,B,C$ on the $\Re$-module $L_n^{(0)}$ are as follows:
\begin{align*}
A u_i^{(0)} &=b_i u_{i-1}^{(0)}+a_i u_i^{(0)}+c_i u_{i+1}^{(0)}
\qquad (0\leq i\leq \floor{\tfrac{n}{2}}),
\\
B u_i^{(0)} &=\theta_i^* u_i^{(0)}
\qquad (0\leq i\leq \floor{\tfrac{n}{2}}),
\\
C u_i^{(0)} &=-b_i u_{i-1}^{(0)}+a_i u_i^{(0)}-c_i u_{i+1}^{(0)}
\qquad (0\leq i\leq \floor{\tfrac{n}{2}}),
\end{align*}
where $u_{-1}$ and $u_{\floor{\tfrac{n}{2}}+1}$ are interpreted as the zero vector and
\begin{align*}
a_i &=\frac{n(n+2)-(n-4i)^2}{32}-\frac{1}{4}
\qquad (0\leq i\leq \floor{\tfrac{n}{2}}),
\\
b_i &=\frac{(n-2i+1)(n-2i+2)}{16}
\qquad (0\leq i\leq \floor{\tfrac{n}{2}}),
\\
c_i &=\frac{(i+1)(2i+1)}{8}
\qquad (0\leq i\leq \floor{\tfrac{n}{2}}),
\\
\theta_i^* &= \frac{(n-4i)^2}{16}-\frac{1}{4}
\qquad (0\leq i\leq \floor{\tfrac{n}{2}}).
\end{align*}
\end{prop}
\begin{proof}
It is routine to verify the proposition by using Lemmas \ref{lem:sharp&ABCD}(i)--(iii) and \ref{lem:Ln0_Ue}.
\end{proof}

\begin{thm}
\label{thm:irr_Ln0}
Suppose that $n\in \N$. Then the following statements hold:
\begin{enumerate}
\item If $n$ is odd then $L_n^{(0)}$ is isomorphic to the irreducible $\Re$-module $R_{\frac{n-1}{2}}(-\frac{1}{4},-\frac{1}{4},-\frac{1}{4})$.

\item If $n=0$ then $L_n^{(0)}$ is isomorphic to the irreducible $\Re$-module $R_0(-\frac{1}{2},-\frac{1}{2},-\frac{1}{2})$.

\item If $n\not=0$ and $n\equiv0\pmod{4}$ then $L_n^{(0)}$ is isomorphic to a direct sum of the irreducible $\Re$-modules $R_{\frac{n}{4}-1}(\frac{n-4}{8},\frac{n}{8},\frac{n-4}{8})$ and
$R_{\frac{n}{4}}(\frac{n-4}{8},\frac{n-4}{8},\frac{n-4}{8})$.

\item If $n\equiv2\pmod{4}$ then $L_n^{(0)}$ is isomorphic to a direct sum of the irreducible $\Re$-modules $R_{\frac{n-2}{4}}(\frac{n-6}{8},\frac{n-2}{8},\frac{n-2}{8})$ and
$R_{\frac{n-2}{4}}(\frac{n-2}{8},\frac{n-2}{8},\frac{n-6}{8})$.
\end{enumerate}
\end{thm}
\begin{proof}
In the proof we employ the notations of Proposition \ref{prop:A,B,C on Ln0}. Note that the scalars $\{b_i\}_{i=0}^{\floor{\frac{n}{2}}}$ and $\{c_i\}_{i=0}^{\floor{\frac{n}{2}}}$ are nonzero.

(i): Suppose that $n$ is odd and set $d=\frac{n-1}{2}\in \N$. The matrices representing $A$ and $C$ with respect to the basis $\{u_i^{(0)}\}_{i=0}^d$ for $L_n^{(0)}$ are  irreducible tridiagonal. The matrix representing $B$ with respect to the basis $\{u_i^{(0)}\}_{i=0}^d$ for $L_n^{(0)}$ is diagonal and the diagonal entries $\{\theta_i^*\}_{i=0}^d$ are mutually distinct. Therefore the $\Re$-module $L_n^{(0)}$ is irreducible.
The traces of $A,B,C$ on $L_n^{(0)}$ are equal to $d+1$ times
$$
-\frac{3}{16}+\frac{d(d+2)}{12}.
$$
By Lemma \ref{lem:Rd(abc)_isomorphism class} the $\Re$-module $L_n^{(0)}$ is isomorphic to $R_d(-\frac{1}{4},-\frac{1}{4},-\frac{1}{4})$.

(ii): Suppose that $n=0$. By Proposition \ref{prop:A,B,C on Ln0} each of $A,B,C$ acts on $L_n^{(0)}$ as the scalar $-\frac{1}{4}$. By Lemma \ref{lem:Rd(abc)_isomorphism class} the $\Re$-module $L_n^{(0)}$ is isomorphic to $R_0(-\frac{1}{2},-\frac{1}{2},-\frac{1}{2})$.

(iii): Suppose that $n\not=0$ and $n\equiv0\pmod{4}$. Set $d=\frac{n}{4}\geq 1$. Let
\begin{align*}
v_i&=u_i^{(0)}-u_{2d-i}^{(0)}
\qquad
(0\leq i\leq d-1),
\\
w_i&=
u_i^{(0)}+u_{2d-i}^{(0)}
\qquad
(0\leq i\leq d).
\end{align*}
Let $V$ and $W$ denote the subspaces of $L_n^{(0)}$ spanned by $\{v_i\}_{i=0}^{d-1}$ and $\{w_i\}_{i=0}^{d}$ respectively.
Since $\{v_i\}_{i=0}^{d-1}$ and $\{w_i\}_{i=0}^d$ form a basis for $L_n^{(0)}$ it follows that
$$
L_n^{(0)}=V\oplus W.
$$
Observe that $a_i=a_{2d-i}$, $b_i=c_{2d-i}$ and $\theta_i^*=\theta_{2d-i}^*$ for all integers $i$ with $0\leq i\leq 2d$. Combined with Proposition \ref{prop:A,B,C on Ln0} this yields that
\begin{align*}
A v_i &=b_i v_{i-1} + a_i v_i + c_i v_{i+1}
\qquad (0\leq i\leq d-1),
\\
B v_i &= \theta_i^* v_i
\qquad (0\leq i\leq d-1),
\\
C v_i &=-b_i v_{i-1} + a_i v_i - c_i v_{i+1}
\qquad (0\leq i\leq d-1),
\end{align*}
where $v_{-1}$ and $v_d$ are interpreted as the zero vector and
\begin{align*}
A w_i &= b_i w_{i-1} + a_i w_i + c_i w_{i+1}
\qquad (0\leq i\leq d-1),
\\
A w_d &= 2 b_d w_{d-1}+a_d w_d,
\\
B w_i &= \theta_i^* w_i
\qquad (0\leq i\leq d),
\\
C w_i &= -b_i w_{i-1} + a_i w_i - c_i w_{i+1}
\qquad (0\leq i\leq d-1),
\\
C w_d &= -2b_d w_{d-1}+a_d w_d.
\end{align*}
where $w_{-1}$ is interpreted as the zero vector. Since the algebra $\Re$ is generated by $A,B,C$ it follows that $V$ and $W$ are two $\Re$-submodules of $L_{n}^{(0)}$.

The matrices representing $A$ and $C$ with respect to the basis $\{v_i\}_{i=0}^{d-1}$ for $V$ are irreducible tridiagonal. The matrix representing $B$ with respect to the basis $\{v_i\}_{i=0}^{d-1}$ for $V$ is diagonal and the diagonal entries $\{\theta_i^*\}_{i=0}^{d-1}$ are mutually distinct. Therefore the $\Re$-module $V$ is irreducible. The trace of $B$ on $V$ is equal to $d$ times
$$
\frac{d(d+2)}{4}+\frac{(d-1)(d+1)}{12}
$$
and the traces of $A$ and $C$ on $V$ are equal to $d$ times
$$
\frac{(d-1)(d+1)}{4}+\frac{(d-1)(d+1)}{12}.
$$
By Lemma \ref{lem:Rd(abc)_isomorphism class} the $\Re$-module $V$
is isomorphic to $R_{d-1}(\frac{d-1}{2},\frac{d}{2},\frac{d-1}{2})$.

The matrices representing $A$ and $C$ with respect to the basis $\{w_i\}_{i=0}^{d}$ for $W$ are irreducible tridiagonal. The matrix representing $B$ with respect to the basis $\{w_i\}_{i=0}^{d}$ for $W$ is diagonal and the diagonal entries $\{\theta_i^*\}_{i=0}^{d}$ are mutually distinct. Therefore the $\Re$-module $W$ is irreducible.
The traces of $A,B,C$ on $W$ are equal to $d+1$ times
\begin{gather*}
\frac{(d-1)(d+1)}{4}+\frac{d(d+2)}{12}.
\end{gather*}
By Lemma \ref{lem:Rd(abc)_isomorphism class} the $\Re$-module $W$ is isomorphic to $R_d(\frac{d-1}{2},\frac{d-1}{2},\frac{d-1}{2})$.

(iv): Suppose that $n\equiv2\pmod{4}$ and set $d=\frac{n-2}{4}\in \N$. Let
\begin{align*}
v_i&=u_i^{(0)}-u_{2d-i+1}^{(0)}
\qquad
(0\leq i\leq d),
\\
w_i&=
u_i^{(0)}+u_{2d-i+1}^{(0)}
\qquad
(0\leq i\leq d).
\end{align*}
Let $V$ and $W$ denote the subspaces of $L_n^{(0)}$ spanned by $\{v_i\}_{i=0}^{d}$ and $\{w_i\}_{i=0}^{d}$ respectively.
Since $\{v_i\}_{i=0}^d$ and $\{w_i\}_{i=0}^d$ form a basis for $L_n^{(0)}$ it follows that
$$
L_n^{(0)}=V\oplus W.
$$
Observe that $a_i=a_{2d-i+1}$, $b_i=c_{2d-i+1}$ and $\theta_i^*=\theta_{2d-i+1}^*$ for all integers $i$ with $0\leq i\leq 2d+1$. Combined with Proposition \ref{prop:A,B,C on Ln0} this yields that
\begin{align*}
A v_i &=
b_i v_{i-1} + a_i v_i + c_i v_{i+1}
\qquad
(0\leq i\leq d-1),
\\
A v_d &= b_d v_{d-1} + (a_d-c_d) v_d,
\\
B v_i &= \theta_i^* v_i
\qquad
(0\leq i\leq d),
\\
C v_i &=-b_i v_{i-1} + a_i v_i - c_i v_{i+1}
\qquad
(0\leq i\leq d-1),
\\
C v_d &=-b_{d} v_{d-1} + (a_{d}+c_d) v_{d},
\end{align*}
where $v_{-1}$ is interpreted as the zero vector and
\begin{align*}
A w_i &=b_i w_{i-1} + a_i w_i + c_i w_{i+1}
\qquad (0\leq i\leq d-1),
\\
A w_d &= b_{d} w_{d-1} + (a_d+c_d) w_d,
\\
B w_i &= \theta_i^* w_i
\qquad (0\leq i\leq d),
\\
C w_i &=  -b_i w_{i-1} + a_i w_i - c_i w_{i+1}
\qquad (0\leq i\leq d-1),
\\
C w_d &= -b_{d} w_{d-1} + (a_d-c_d) w_d,
\end{align*}
where $w_{-1}$ is interpreted as the zero vector. Since the algebra $\Re$ is generated by $A,B,C$ it follows that $V$ and $W$ are two $\Re$-submodules of $L_n^{(0)}$.

The matrices representing $A$ and $C$ with respect to the basis $\{v_i\}_{i=0}^d$ for $V$ are irreducible tridiagonal. The matrix representing $B$ with respect to the basis $\{v_i\}_{i=0}^d$ for $V$ is diagonal and the diagonal entries $\{\theta_i^*\}_{i=0}^d$ are mutually distinct.
Therefore the $\Re$-module $V$ is irreducible.
The traces of $A,B,C$ on $V$ are equal to $d+1$ times
\begin{gather*}
\frac{(d-1)(d+1)}{4}+\frac{d(d+2)}{12},
\qquad
\frac{d(d+2)}{4}+\frac{d(d+2)}{12},
\qquad
\frac{d(d+2)}{4}+\frac{d(d+2)}{12},
\end{gather*}
respectively.
By Lemma \ref{lem:Rd(abc)_isomorphism class} the $\Re$-module $V$
is isomorphic to $R_{d}(\frac{d-1}{2},\frac{d}{2},\frac{d}{2})$.

The matrices representing $A$ and $C$ with respect to the basis $\{w_i\}_{i=0}^d$ for $W$ are irreducible tridiagonal. The matrix representing $B$ with respect to the basis $\{w_i\}_{i=0}^d$ for $W$ is diagonal and the diagonal entries $\{\theta_i^*\}_{i=0}^d$ are mutually distinct.
Therefore the $\Re$-module $W$ is irreducible.
The traces of $A,B,C$ on $W$ are equal to $d+1$ times
\begin{gather*}
\frac{d(d+2)}{4}+\frac{d(d+2)}{12},
\qquad
\frac{d(d+2)}{4}+\frac{d(d+2)}{12},
\qquad
\frac{(d-1)(d+1)}{4}+\frac{d(d+2)}{12},
\end{gather*}
respectively.
By Lemma \ref{lem:Rd(abc)_isomorphism class} the $\Re$-module $W$ is isomorphic to $R_{d}(\frac{d}{2},\frac{d}{2},\frac{d-1}{2})$.
\end{proof}

\begin{lem}
[Lemma 5.7, \cite{halved:2023}]
\label{lem:Ln1_Ue}
For any integer $n\geq 1$ the $\Ue$-module $L_n^{(1)}$ satisfies the following conditions:
\begin{enumerate}
\item   There exists a basis $\{u_i^{(1)}\}_{i=0}^{\floor{\frac{n-1}{2}}}$ for $L_n^{(1)}$ such that
        \begin{align*}
         E^2 u_{i}^{(1)} &=(n-2i)(n-2i+1) u_{i-1}^{(1)}
        \quad
        (1\leq i\leq \textstyle{\floor{\frac{n-1}{2}}}),
        \qquad
        E^2 u_0^{(1)}=0,
        \\
         F^2 u_{i}^{(1)} &=(2i+2)(2i+3) u_{i+1}^{(1)}
        \quad
        (0\leq i\leq \textstyle\floor{\frac{n-1}{2}}-1),
        \qquad
        F^2 u_{\floor{\frac{n-1}{2}}}^{(1)}=0,
        \\
         H u_{i}^{(1)} &=(n-4i-2) u_{i}^{(1)}
        \quad
        (0\leq i\leq \textstyle\floor{\frac{n-1}{2}}).
\end{align*}
\item   The element $\Lambda$ acts on $L_n^{(1)}$ as scalar multiplication by $\frac{n(n+2)}{2}$.
\end{enumerate}
\end{lem}

The notation $\{u_i^{(1)}\}_{i=0}^{\floor{\frac{n-1}{2}}}$ will always stand for the basis for the $\U_e$-module $L_n^{(1)}$ ($n\geq 1$) given in Lemma \ref{lem:Ln1_Ue}.
By Proposition \ref{prop:sharp&Ue} the $\U_e$-module $L_n^{(1)}$ is an $\Re$-module by pulling back via $\sharp$.

\begin{prop}\label{prop:A,B,C on Ln1}
Suppose that $n\geq 1$ is an integer. Then the actions of $A,B,C$ on the $\Re$-module $L_n^{(1)}$ are as follows:
\begin{align*}
A u_i^{(1)} &= b_i u_{i-1}^{(1)}+a_i u_i^{(1)}+c_i u_{i+1}^{(1)}
\qquad
(0\leq i\leq \floor{\tfrac{n-1}{2}}),
\\
B u_i^{(1)} &=\theta_i^* u_i^{(1)}
\qquad
(0\leq i\leq \floor{\tfrac{n-1}{2}}),
\\
C u_i^{(1)} &=
-b_i u_{i-1}^{(1)}+a_i u_i^{(1)}-c_i u_{i+1}^{(1)}
\qquad
(0\leq i\leq \floor{\tfrac{n-1}{2}}),
\end{align*}
where $u_{-1}$ and $u_{\floor{\tfrac{n+1}{2}}}$ are interpreted as the zero vector and
\begin{align*}
a_i &=\frac{n(n+2)-(n-4i-2)^2}{32}-\frac{1}{4}
\qquad (0\leq i\leq \floor{\tfrac{n-1}{2}}),
\\
b_i &=\frac{(n-2i)(n-2i+1)}{16}
\qquad (0\leq i\leq \floor{\tfrac{n-1}{2}}),
\\
c_i &=\frac{(i+1)(2i+3)}{8}
\qquad (0\leq i\leq \floor{\tfrac{n-1}{2}}),
\\
\theta_i^* &= \frac{(n-4i-2)^2}{16}-\frac{1}{4}
\qquad (0\leq i\leq \floor{\tfrac{n-1}{2}}).
\end{align*}
\end{prop}
\begin{proof}
It is routine to verify the proposition by using Lemmas \ref{lem:sharp&ABCD}(i)--(iii) and \ref{lem:Ln1_Ue}.
\end{proof}

\begin{thm}
\label{thm:irr_Ln1}
Suppose that $n\geq 1$ is an integer. Then the following statements hold:
\begin{enumerate}
\item If $n$ is odd then $L_n^{(1)}$ is isomorphic to the irreducible $\Re$-module $
R_{\frac{n-1}{2}}(
-\frac{1}{4},-\frac{1}{4},-\frac{1}{4})$.

\item If $n=2$ then $L_n^{(1)}$ is isomorphic to the irreducible $\Re$-module $R_0(0,-\frac{1}{2},0)$.

\item If $n\not=2$ and $n\equiv2\pmod{4}$ then $L_n^{(1)}$ is isomorphic to a direct sum of the irreducible $\Re$-modules $R_{\frac{n-6}{4}}(\frac{n-2}{8},\frac{n-2}{8},\frac{n-2}{8})$ and
$R_{\frac{n-2}{4}}(\frac{n-2}{8},\frac{n-6}{8},\frac{n-2}{8})$.

\item If $n\equiv0\pmod{4}$ then $L_n^{(1)}$ is isomorphic to a direct sum of the irreducible $\Re$-modules $R_{\frac{n}{4}-1}(\frac{n}{8},\frac{n-4}{8},\frac{n-4}{8})$ and $R_{\frac{n}{4}-1}(\frac{n-4}{8},\frac{n-4}{8},\frac{n}{8})$.
\end{enumerate}
\end{thm}
\begin{proof}
In the proof we employ the notations of Proposition \ref{prop:A,B,C on Ln1}. Note that the scalars $\{b_i\}_{i=0}^{\floor{\frac{n-1}{2}}}$ and $\{c_i\}_{i=0}^{\floor{\frac{n-1}{2}}}$ are nonzero.

(i): Suppose that $n$ is odd and set $d=\frac{n-1}{2}\in \N$.
The matrices representing $A$ and $C$ with respect to the basis $\{u_i^{(1)}\}_{i=0}^d$ for $L_n^{(1)}$ are  irreducible tridiagonal. The matrix representing $B$ with respect to the basis $\{u_i^{(1)}\}_{i=0}^d$ for $L_n^{(1)}$ is diagonal and the diagonal entries $\{\theta_i^*\}_{i=0}^d$ are mutually distinct. Therefore the $\Re$-module $L_n^{(1)}$ is irreducible.
The traces of $A,B,C$ on $L_n^{(1)}$ are equal to $d+1$ times
$$
-\frac{3}{16}+\frac{d(d+2)}{12}.
$$
By Lemma \ref{lem:Rd(abc)_isomorphism class} the $\Re$-module $L_n^{(1)}$ is isomorphic to $R_d(-\frac{1}{4},-\frac{1}{4},-\frac{1}{4})$.

(ii): Suppose that $n=2$. By Proposition \ref{prop:A,B,C on Ln1} the elements $A,B,C$ act on $L_n^{(1)}$ as the scalars $0,-\frac{1}{4},0$ respectively. By Lemma \ref{lem:Rd(abc)_isomorphism class} the $\Re$-module $L_n^{(1)}$ is isomorphic to $R_0(0,-\frac{1}{2},0)$.

(iii): Suppose that $n\not=2$ and $n\equiv2\pmod{4}$. Set $d=\frac{n-2}{4}\geq 1$. Let
\begin{align*}
v_i&=u_i^{(1)}-u_{2d-i}^{(1)}
\qquad (0\leq i\leq d-1),
\\
w_i&=
u_i^{(1)}+u_{2d-i}^{(1)}
\qquad (0\leq i\leq d).
\end{align*}
Let $V$ and $W$ denote the subspaces of $L_{n}^{(1)}$ spanned by $\{v_i\}_{i=0}^{d-1}$ and $\{w_i\}_{i=0}^{d}$ respectively.
Since $\{v_i\}_{i=0}^{d-1}$ and $\{w_i\}_{i=0}^d$ form a basis for $L_{n}^{(1)}$ it follows that
$$
L_{n}^{(1)}=V\oplus W.
$$
Observe that $a_i=a_{2d-i}$, $b_i=c_{2d-i}$ and $\theta_i^*=\theta_{2d-i}^*$ for all integers $i$ with $0\leq i\leq 2d$. Combined with Proposition \ref{prop:A,B,C on Ln1} this yields that
\begin{align*}
A v_i &=
b_i v_{i-1} + a_i v_i + c_i v_{i+1}
\qquad
(0\leq i\leq d-1),
\\
B v_i &= \theta_i^* v_i
\qquad
(0\leq i\leq d-1),
\\
C v_i &=
-b_i v_{i-1} + a_i v_i - c_i v_{i+1}
\qquad
(0\leq i\leq d-1),
\end{align*}
where $v_{-1}$ and $v_d$ are interpreted as the zero vector and
\begin{align*}
A w_i &=
b_i w_{i-1} + a_i w_i + c_i w_{i+1}
\qquad
(0\leq i\leq d-1),
\\
A w_d &= 2 b_d w_{d-1} + a_d w_d,
\\
B w_i &= \theta_i^* w_i
\qquad
(0\leq i\leq d),
\\
C w_i &=
-b_i w_{i-1} + a_i w_i - c_i w_{i+1}
\qquad
(0\leq i\leq d-1),
\\
C w_d &= -2 b_{d} w_{d-1} + a_d w_d,
\end{align*}
where $w_{-1}$ is interpreted as the zero vector. Since the algebra $\Re$ is generated by $A,B,C$ it follows that $V$ and $W$ are two $\Re$-submodules of $L_n^{(1)}$.

The matrices representing $A$ and $C$ with respect to the basis $\{v_i\}_{i=0}^{d-1}$ for $V$ are irreducible tridiagonal. The matrix representing $B$ with respect to the basis $\{v_i\}_{i=0}^{d-1}$ for $V$ is diagonal and the diagonal entries $\{\theta_i^*\}_{i=0}^{d-1}$ are mutually distinct. Therefore the $\Re$-module $V$ is irreducible. The traces of $A,B,C$ on $V$ are equal to $d$ times
\begin{gather*}
\frac{d(d+2)}{4}+\frac{(d-1)(d+1)}{12}.
\end{gather*}
By Lemma \ref{lem:Rd(abc)_isomorphism class} the $\Re$-module $V$
is isomorphic to $R_{d-1}(\frac{d}{2},\frac{d}{2},\frac{d}{2})$.

The matrices representing $A$ and $C$ with respect to the basis $\{w_i\}_{i=0}^{d}$ for $W$ are irreducible tridiagonal. The matrix representing $B$ with respect to the basis $\{w_i\}_{i=0}^{d}$ for $W$ is diagonal and the diagonal entries $\{\theta_i^*\}_{i=0}^{d}$ are mutually distinct. Therefore the $\Re$-module $W$ is irreducible.
The traces of $A,B,C$ on $W$ are equal to $d+1$ times
\begin{gather*}
\frac{d(d+2)}{4}+\frac{d(d+2)}{12},
\qquad
\frac{(d-1)(d+1)}{4}+\frac{d(d+2)}{12},
\qquad
\frac{d(d+2)}{4}+\frac{d(d+2)}{12},
\end{gather*}
respectively.
By Lemma \ref{lem:Rd(abc)_isomorphism class} the $\Re$-module $W$ is isomorphic to $R_{d}(\frac{d}{2},\frac{d-1}{2},\frac{d}{2})$.

(iv): Suppose that $n\equiv0\pmod{4}$ and set $d=\frac{n}{4}-1\in \N$. Let
\begin{align*}
v_i&=u_i^{(1)}-u_{2d-i+1}^{(1)}
\qquad
(0\leq i\leq d),
\\
w_i&=
u_i^{(1)}+u_{2d-i+1}^{(1)}
\qquad
(0\leq i\leq d).
\end{align*}
Let $V$ and $W$ denote the subspaces of $L_{n}^{(1)}$ spanned by $\{v_i\}_{i=0}^{d}$ and $\{w_i\}_{i=0}^{d}$ respectively.
Since $\{v_i\}_{i=0}^d$ and $\{w_i\}_{i=0}^d$ form a basis for $L_{n}^{(1)}$ it follows that
$$
L_{n}^{(1)}=V\oplus W.
$$
Observe that $a_i=a_{2d-i+1}$, $b_i=c_{2d-i+1}$ and $\theta_i^*=\theta_{2d-i+1}^*$ for all integers $i$ with $0\leq i\leq 2d+1$. Combined with Proposition \ref{prop:A,B,C on Ln1} this yields that
\begin{align*}
A v_i &=
b_i v_{i-1} + a_i v_i + c_i v_{i+1}
\qquad
(0\leq i\leq d-1),
\\
A v_d &= b_d v_{d-1} + (a_d-c_d) v_d,
\\
B v_i &= \theta_i^* v_i
\qquad
(0\leq i\leq d),
\\
C v_i &=
-b_i v_{i-1} + a_i v_i - c_i v_{i+1}
\qquad
(0\leq i\leq d-1),
\\
C v_d &= -b_{d} v_{d-1} + (a_{d}+c_d) v_d,
\end{align*}
where $v_{-1}$ is interpreted as the zero vector and
\begin{align*}
A w_i &=
b_i w_{i-1} + a_i w_i + c_i w_{i+1}
\qquad
(0\leq i\leq d-1),
\\
A w_d &= b_{d} w_{d-1} + (a_d+c_d) w_d,
\\
B w_i &= \theta_i^* w_i
\qquad
(0\leq i\leq d),
\\
C w_i &=
-b_i w_{i-1} + a_i w_i - c_i w_{i+1}
\qquad
(0\leq i\leq d-1),
\\
C w_d &= -b_{d} w_{d-1} + (a_d-c_d) w_d,
\end{align*}
where $w_{-1}$ is interpreted as the zero vector.
Since the algebra $\Re$ is generated by $A,B,C$ it follows that $V$ and $W$ are two $\Re$-submodules of $L_{n}^{(1)}$.

The matrices representing $A$ and $C$ with respect to the basis $\{v_i\}_{i=0}^d$ for $V$ are irreducible tridiagonal. The matrix representing $B$ with respect to the basis $\{v_i\}_{i=0}^d$ for $V$ is diagonal and the diagonal entries $\{\theta_i^*\}_{i=0}^d$ are mutually distinct.
Therefore the $\Re$-module $V$ is irreducible.
The traces of $A,B,C$ on $V$ are equal to $d+1$ times
\begin{gather*}
\frac{d(d+2)}{4}+\frac{d(d+2)}{12},
\qquad
\frac{d(d+2)}{4}+\frac{d(d+2)}{12},
\qquad
\frac{(d+1)(d+3)}{4}+\frac{d(d+2)}{12},
\end{gather*}
respectively. By Lemma \ref{lem:Rd(abc)_isomorphism class} the $\Re$-module $V$
is isomorphic to $R_{d}(\frac{d}{2},\frac{d}{2},\frac{d+1}{2})$.

The matrices representing $A$ and $C$ with respect to the basis $\{w_i\}_{i=0}^d$ for $W$ are irreducible tridiagonal. The matrix representing $B$ with respect to the basis $\{w_i\}_{i=0}^d$ for $W$ is diagonal and the diagonal entries $\{\theta_i^*\}_{i=0}^d$ are mutually distinct.
Therefore the $\Re$-module $W$ is irreducible.
The traces of $A,B,C$ on $W$ are equal to $d+1$ times
\begin{gather*}
\frac{(d+1)(d+3)}{4}+\frac{d(d+2)}{12},
\qquad
\frac{d(d+2)}{4}+\frac{d(d+2)}{12},
\qquad
\frac{d(d+2)}{4}+\frac{d(d+2)}{12},
\end{gather*}
respectively.
By Lemma \ref{lem:Rd(abc)_isomorphism class} the $\Re$-module $W$ is isomorphic to $R_{d}(\frac{d+1}{2},\frac{d}{2},\frac{d}{2})$.
\end{proof}

\noindent {\it Proof of Theorem \ref{thm:R-module decomposition}.}
(i): Based on the first two paragraphs of this section, it is enough to study the $\U_e$-module $L_n^{(0)}$ ($n\in \N$) and the $\U_e$-module $L_n^{(1)}$ ($n\geq 1$) from the perspective of $\Re$-modules. Therefore (i) is immediate from Theorems \ref{thm:irr_Ln0} and \ref{thm:irr_Ln1}.

(ii): Applying Theorem~\ref{thm:LT} it is routine to verify that $A,B,C$ act as a Leonard triple on every finite-dimensional irreducible $\Re$-module listed in Theorems \ref{thm:irr_Ln0} and \ref{thm:irr_Ln1}. The statement (ii) follows.
\hfill $\square$

\begin{prop}
\label{prop:irr_Ln0&Ln1}
The finite-dimensional irreducible $\Re$-modules mentioned in Theorems {\rm \ref{thm:irr_Ln0}(i)--(iv)} and {\rm \ref{thm:irr_Ln1}(ii)--(iv)}, namely
\begin{align*}
&\textstyle R_{\frac{n-1}{2}}(-\frac{1}{4},-\frac{1}{4},-\frac{1}{4})
\qquad
\hbox{for all odd integers $n\geq 1$};
\\
&\textstyle
R_{\frac{n-2}{4}} (\frac{n-2}{8}, \frac{n-2}{8}, \frac{n-6}{8})
\qquad
\hbox{for all integers $n\geq 2$ with $n\equiv 2 \!\!\!\pmod{4}$};
\\
&\textstyle
R_{\frac{n-2}{4}} (\frac{n-2}{8}, \frac{n-6}{8}, \frac{n-2}{8})
\qquad
\hbox{for all integers $n\geq 2$ with $n\equiv 2 \!\!\!\pmod{4}$};
\\
&\textstyle
R_{\frac{n-2}{4}} (\frac{n-6}{8}, \frac{n-2}{8}, \frac{n-2}{8})
\qquad
\hbox{for all integers $n\geq 2$ with $n\equiv 2 \!\!\!\pmod{4}$};
\\
&\textstyle
R_{\frac{n-6}{4}}(\frac{n-2}{8}, \frac{n-2}{8}, \frac{n-2}{8})
\qquad
\hbox{for all integers $n\geq 6$ with $n\equiv 2 \!\!\!\pmod{4}$};
\\
&\textstyle
R_{\frac{n}{4}-1}(\frac{n-4}{8}, \frac{n-4}{8}, \frac{n}{8})
\qquad
\hbox{for all integers $n\geq 4$ with $n\equiv 0 \!\!\! \pmod{4}$};
\\
&\textstyle
R_{\frac{n}{4}-1}(\frac{n-4}{8}, \frac{n}{8}, \frac{n-4}{8})
\qquad
\hbox{for all integers $n\geq 4$ with $n\equiv 0 \!\!\! \pmod{4}$};
\\
&\textstyle
R_{\frac{n}{4}-1}(\frac{n}{8}, \frac{n-4}{8}, \frac{n-4}{8})
\qquad
\hbox{for all integers $n\geq 4$ with $n\equiv 0 \!\!\!\pmod{4}$};
\\
&\textstyle
R_{\frac{n}{4}}(\frac{n-4}{8}, \frac{n-4}{8}, \frac{n-4}{8})
\qquad
\hbox{for all integers $n\geq 0$ with $n\equiv 0 \!\!\!\pmod{4}$},
\end{align*}
are mutually non-isomorphic.
\end{prop}
\begin{proof}
It is straightforward to verify the proposition by applying Lemma \ref{lem:Rd(abc)_isomorphism class}.
\end{proof}

We finish this section with an improvement of Proposition \ref{prop:sharp&Ue}.

\begin{thm}
${\rm Im}
\, \sharp$ is a proper subalgebra of $\U_e$.
\end{thm}
\begin{proof}
Suppose that ${\rm Im}\, \sharp=\U_e$. Then every irreducible $\U_e$-module is an irreducible $\Re$-module, a contradiction to Theorems \ref{thm:irr_Ln0}(iii), (iv) and \ref{thm:irr_Ln1}(iii), (iv).
Suppose that ${\rm Im}\, \sharp$ is the subalgebra of $\U_e$ generated by the unit.
Then every irreducible ${\rm Im}\, \sharp$-module is of dimension one, a contradiction to Theorems \ref{thm:irr_Ln0}(i) and \ref{thm:irr_Ln1}(i).
\end{proof}

\section{Proof for Theorem \ref{thm:LP_hypercube+johnson}}
\label{sec:Proof_for_thm_LP_hypercube+johnson}

Fix an integer $D\geq 2$. Recall from \S\ref{sec:intro} that the notation $X$ denotes the power set of a $D$-element set and the notation $H(D,2)$ represents the $D$-dimensional hypercube.
For any integer $i$ with $0\leq i\leq D$ let $R_i$ denote the $i^{\rm\, th}$ distance relation in $H(D,2)$.

\begin{lem}
[Theorem 13.2, \cite{hypercube2002}]
\label{lem:U_module_hypercube}
There exists a unique $\U$-module $\C^{X}$ such that
  \begin{align*}
    Ex &= \sum_{\substack{|y|<|x|\\ (x,y)\in R_1}} y  \qquad \hbox{for all $x\in X$},  \\
    Fx &= \sum_{\substack{|y|>|x|\\ (x,y)\in R_1}}  y  \qquad \hbox{for all $x\in X$}, \\
    Hx &= (D-2|x|)x \qquad \hbox{for all $x\in X$}.
  \end{align*}
\end{lem}

Recall the algebra homomorphism $\sharp:\Re\to \U$ from Theorem \ref{thm:RetoU}. An $\Re$-module $\C^{X}$ is obtained by pulling back the $\U$-module $\C^{X}$ via $\sharp$.
Recall from Proposition \ref{prop:sharp&Ue} that ${\rm Im}\,\sharp$ is contained in the even subalgebra $\U_e$ of $\U$. By Lemma \ref{lem:Ue_basis} the algebra $\U_e$ is generated by $E^2,F^2,H,\Lambda$.
The $\U_e$-module $\C^{X}$ is as follows:

\begin{lem}
\label{lem:Ue_module_hypercube}
There exists a unique $\U_e$-module $\C^{X}$ such that
  \begin{align}
    E^2 x &=  2\sum_{\substack{|y|<|x|\\ (x,y)\in R_2}} y
    \qquad
    \hbox{for all $x\in X$},
    \label{E2_hypercube}
    \\
    F^2 x &= 2 \sum_{\substack{|y|>|x|\\ (x,y)\in R_2}} y
    \qquad
    \hbox{for all $x\in X$},
    \label{F2_hypercube}
    \\
    Hx &= (D-2|x|)x
    \qquad \hbox{for all $x\in X$},
    \label{H_hypercube}
    \\
    \Lambda x &=
     \left(D+\frac{(D-2|x|)^2}{2}\right) x
     +2\sum_{\substack{|y|=|x|\\ (x,y)\in R_2}} y
     \qquad \hbox{for all $x\in X$}.
     \label{Lambda_hypercube}
  \end{align}
\end{lem}
\begin{proof}
The equation (\ref{H_hypercube}) is immediate from Lemma \ref{lem:U_module_hypercube}. The equation (\ref{Lambda_hypercube}) is immediate from \cite[Lemma 5.2]{Huang:CG&Johnson}.
Let $x\in X$ be given. For any $y\in X$ with $|y|<|x|$ and $(x,y)\in R_2$, there are exactly two paths of length two from $x$ to $y$ in $H(D,2)$. Hence (\ref{E2_hypercube}) follows. For any $y\in X$ with $|y|>|x|$ and $(x,y)\in R_2$, there are exactly two paths of length two from $x$ to $y$ in $H(D,2)$. Hence (\ref{F2_hypercube}) follows.
\end{proof}

It is now easy to evaluate the actions of the generators $A,B,C$ of $\Re$ on the $\Re$-module $\C^{X}$ by applying Lemma \ref{lem:Ue_module_hypercube} to Lemma \ref{lem:sharp&ABCD}(i)--(iii):

\begin{lem}\label{lem:R_module_hypercube}
There exists a unique $\Re$-module $\C^{X}$ such that
\begin{align*}
  Ax &=
  \left(\frac{D}{16}-\frac{1}{4}\right) x
  +\frac{1}{8}\sum_{(x,y)\in R_2}y
  \qquad \hbox{for all $x\in X$},
  \\
  Bx &=
  \left(\frac{(D-2|x|)^2}{16}-\frac{1}{4}\right) x
  \qquad \hbox{for all $x\in X$},
  \\
  Cx &=
  \left(\frac{D}{16}-\frac{1}{4}\right) x
  +\frac{1}{8}\sum_{\substack{|y|=|x| \\ (x,y)\in R_2}} y
  -\frac{1}{8}\sum_{\substack{|y|\not=|x| \\ (x,y)\in R_2}} y
  \qquad \hbox{for all $x\in X$}.
\end{align*}
\end{lem}

\noindent {\it Proof of Theorem \ref{thm:LP_hypercube+johnson}.}
In the proof the notation $\phi$ is used to indicate the algebra homomorphism $\Re\to \R$ described in Theorem \ref{thm:LP_hypercube+johnson}(i). The existence and uniqueness of $\phi$ are proved as follows:

(i): Recall the linear maps $\A_2^J, \A_2^{\bar J}, \A_2^*$ on $\C^{X}$ from (\ref{eq:A2+})--(\ref{eq:dual_A_hypercube}).
By Lemma \ref{lem:R_module_hypercube} the actions of $A,B,C$ on the $\Re$-module $\C^{X}$ are identical to the linear maps given in the right-hand sides of (\ref{A:R->R})--(\ref{C:R->R}) respectively. The existence of $\phi$ follows. Since the algebra $\Re$ is generated by $A,B,C$ the uniqueness of $\phi$ follows. The linear maps $\A_2^J$, $\A_2^{\bar J}$, $\A_2^*$ are identical to the actions of
\begin{equation*}
2-\frac{D}{2}+4(A+C),
\qquad
4(A-C),
\qquad
2-\frac{D}{2}+8B
\end{equation*}
on the $\Re$-module $\C^{X}$ respectively. Since the algebra $\R$ is generated by $\A_2^J, \A_2^{\bar J}, \A_2^*$ the algebra homomorphism $\phi$ is onto. The statement (i) follows.

(ii): Recall from the paragraph below Theorem \ref{thm:R-module decomposition} that the algebra $\R$ is a finite-dimensional semisimple algebra and every irreducible $\R$-module is contained in the standard $\R$-module $\C^{X}$ up to isomorphism. By Theorem \ref{thm:R-module decomposition} the $\Re$-module $\C^{X}$ is completely reducible and $A,B,C$ act on every irreducible $\Re$-submodule of $\C^{X}$ as a Leonard triple. Combined with the surjectivity of $\phi$ the statement (ii) follows.
\hfill $\square$

\section{The algebraic properties of $\R$}\label{sec:R_structure}

Recall the $(n+1)$-dimensional irreducible $\U$-module $L_n$ ($n\in \N$) from the first paragraph of \S\ref{sec:Remodule_Ln}.
Note that every finite-dimensional $\U$-module is completely reducible.

\begin{lem}
\label{lem:dec_U-CX}
\begin{enumerate}
\item Suppose that $D$ is odd. Then the pairwise non-isomorphic irreducible $\U$-submodules of $\C^{X}$ are
$$
L_{2k+1}
\qquad \textstyle (0\leq k\leq \frac{D-1}{2}).
$$
\item Suppose that $D$ is even. Then the pairwise non-isomorphic irreducible $\U$-submodules of $\C^{X}$ are
$$
L_{2k}
\qquad \textstyle (0\leq k\leq \frac{D}{2}).
$$
\end{enumerate}
\end{lem}
\begin{proof}
Immediate from \cite[Theorem 10.2]{hypercube2002} or \cite[Theorem 5.3]{Huang:CG&Johnson}.
\end{proof}

Recall the $(\floor{\frac{n}{2}}+1)$-dimensional irreducible $\U_e$-module $L_n^{(0)}$ ($n\in \N$) and
the $\floor{\frac{n+1}{2}}$-dimensional irreducible $\U_e$-module $L_n^{(1)}$ ($n\geq 1$) from the second paragraph of \S\ref{sec:Remodule_Ln}.

\begin{lem}
\label{lem:dec_Ue-CX}
\begin{enumerate}
\item Suppose that $D$ is odd. Then the pairwise non-isomorphic irreducible $\U_e$-submodules  of $\C^X$ are
\begin{align*}
L_{2k+1}^{(0)}
\qquad \textstyle (0\leq k\leq \frac{D-1}{2});
\qquad
L_{2k+1}^{(1)}
\qquad \textstyle (0\leq k\leq \frac{D-1}{2}).
\end{align*}
\item Suppose that $D$ is even. Then the pairwise non-isomorphic irreducible $\U_e$-submodules of $\C^X$ are
\begin{align*}
L_{2k}^{(0)}
\qquad \textstyle (0\leq k\leq \frac{D}{2});
\qquad
L_{2k}^{(1)}
\qquad \textstyle (1\leq k\leq \frac{D}{2}).
\end{align*}
\end{enumerate}
\end{lem}
\begin{proof}
Apply (\ref{eq:Ln_decomposition_as_U_e}) to Lemma \ref{lem:dec_U-CX}.
\end{proof}

Recall the finite-dimensional irreducible $\Re$-modules from \S\ref{sec:Remodule}.

\begin{lem}
\label{lem:dec_R-CX}
\begin{enumerate}
\item Suppose that $D$ is odd. Then the pairwise non-isomorphic irreducible $\Re$-submodules of $\C^X$ are
$$
\textstyle
R_k(-\frac{1}{4},-\frac{1}{4},-\frac{1}{4})
\qquad  (0\leq k\leq \frac{D-1}{2}).
$$
\item Suppose that $D\equiv2 \pmod{4}$. Then the pairwise non-isomorphic irreducible $\Re$-submodules of $\C^X$ are
\begin{alignat*}{2}
&\textstyle R_k(\frac{k}{2},\frac{k+1}{2},\frac{k}{2})
\qquad
(0\leq k\leq \frac{D-6}{4});
\qquad
&&\textstyle R_{k}(\frac{k-1}{2},\frac{k-1}{2},\frac{k-1}{2})
\qquad
(0\leq k\leq \frac{D-2}{4});
\\
&\textstyle R_{k}(\frac{k-1}{2},\frac{k}{2},\frac{k}{2})
\qquad
(0\leq k\leq \frac{D-2}{4});
\qquad
&&\textstyle R_{k}(\frac{k}{2},\frac{k}{2},\frac{k-1}{2})
\qquad
(0\leq k\leq \frac{D-2}{4});
\\
&\textstyle R_k(\frac{k+1}{2},\frac{k+1}{2},\frac{k+1}{2})
\qquad
(0\leq k\leq \frac{D-6}{4});
\qquad
&&\textstyle R_k(\frac{k}{2},\frac{k-1}{2},\frac{k}{2})
\qquad
(0\leq k\leq \frac{D-2}{4});
\\
&\textstyle R_k(\frac{k+1}{2},\frac{k}{2},\frac{k}{2})
\qquad
(0\leq k\leq \frac{D-6}{4});
\qquad
&&\textstyle R_k(\frac{k}{2},\frac{k}{2},\frac{k+1}{2})
\qquad
(0\leq k\leq \frac{D-6}{4}).
\end{alignat*}
\item Suppose that $D\equiv0 \pmod{4}$. Then the pairwise non-isomorphic irreducible $\Re$-submodules of $\C^X$ are
\begin{alignat*}{2}
&\textstyle R_k(\frac{k}{2},\frac{k+1}{2},\frac{k}{2})
\qquad
(0\leq k\leq \frac{D}{4}-1);
\qquad
&&\textstyle R_{k}(\frac{k-1}{2},\frac{k-1}{2},\frac{k-1}{2})
\qquad
(0\leq k\leq \frac{D}{4});
\\
&\textstyle R_{k}(\frac{k-1}{2},\frac{k}{2},\frac{k}{2})
\qquad
(0\leq k\leq \frac{D}{4}-1);
\qquad
&&\textstyle R_{k}(\frac{k}{2},\frac{k}{2},\frac{k-1}{2})
\qquad
(0\leq k\leq \frac{D}{4}-1);
\\
&\textstyle R_k(\frac{k+1}{2},\frac{k+1}{2},\frac{k+1}{2})
\qquad
(0\leq k\leq \frac{D}{4}-2);
\qquad
&&\textstyle R_k(\frac{k}{2},\frac{k-1}{2},\frac{k}{2})
\qquad
(0\leq k\leq \frac{D}{4}-1);
\\
&\textstyle R_k(\frac{k+1}{2},\frac{k}{2},\frac{k}{2})
\qquad
(0\leq k\leq \frac{D}{4}-1);
\qquad
&&\textstyle R_k(\frac{k}{2},\frac{k}{2},\frac{k+1}{2})
\qquad
(0\leq k\leq \frac{D}{4}-1).
\end{alignat*}
\end{enumerate}
\end{lem}
\begin{proof}
Apply Theorems \ref{thm:irr_Ln0} and \ref{thm:irr_Ln1} along with Proposition \ref{prop:irr_Ln0&Ln1} to Lemma \ref{lem:dec_Ue-CX}.
\end{proof}

Recall the algebra $\R$ from the paragraph below Theorem \ref{thm:R-module decomposition}.

\begin{thm}\label{thm:R_structure}
\begin{enumerate}
\item The algebra $\R$ is isomorphic to
\begin{align*}
\left\{
\begin{array}{ll}
\bigoplus\limits_{k=1}^{\frac{D+1}{2}}\End (\C^k)
\qquad
\hbox{if $D$ is odd},
\\
4\cdot \End(\C^{\frac{D+2}{4}})\oplus
8\cdot \bigoplus\limits_{k=1}^{\frac{D-2}{4}} \End(\C^k)
\qquad
\hbox{if $D\equiv 2 \!\!\! \pmod{4}$},
\\
\End(\C^{\frac{D}{4}+1})
\oplus 7\cdot\End(\C^{\frac{D}{4}})
\oplus
8\cdot \bigoplus\limits_{k=1}^{\frac{D}{4}-1}\End(\C^k)
\qquad
\hbox{if $D\equiv 0 \!\!\! \pmod{4}$}.
\end{array}
\right.
\end{align*}
\item The dimension of $\R$ is equal to
$
{\floor{\frac{D}{2}}+3\choose 3}+
{\ceil{\frac{D}{2}}+1\choose 3}.
$
\end{enumerate}
\end{thm}
\begin{proof}
(i): By Theorem \ref{thm:LP_hypercube+johnson}(i) and since $\R$ is a semisimple subalgebra of $\End(\C^X)$, the algebra $\R$ is isomorphic to
  $$
  \bigoplus\End(\C^{\dim V})
  $$
where the direct sum is over all non-isomorphic irreducible $\Re$-submodules $V$ of $\C^{X}$.
Therefore (i) is immediate from Lemma \ref{lem:dec_R-CX}.

(ii): By Theorem \ref{thm:R_structure}(i) the dimension of $\R$ is equal to
\begin{gather*}
\left\{
\begin{array}{ll}
\ds \sum_{k=1}^{\frac{D+1}{2}} k^2
\qquad
&\hbox{if $D$ is odd},
\\
\ds \frac{(D+2)^2}{4}
+8\sum_{k=1}^{\frac{D-2}{4}} k^2
\qquad
&\hbox{if $D\equiv 2 \!\!\! \pmod{4}$},
\\
\ds \frac{(D+4)^2}{16}
+\frac{7 D^2}{16}
+
8\sum_{k=1}^{\frac{D}{4}-1} k^2
\qquad
&\hbox{if $D\equiv 0  \!\!\! \pmod{4}$}.
\end{array}
\right.
\end{gather*}
For each case a direct calculation shows that (ii) holds.
\end{proof}

Let $X_e$ denote the equivalence class of the empty set $\emptyset$ in $X$ with respect to the equivalence relation generated by $R_2$. More explicitly
$$
X_e=\{x\in X\,|\, \hbox{$|x|$ is even}\}.
$$
The {\it halved cube} $\frac{1}{2}H(D,2)$ is a finite simple connected graph with vertex set $X_e$ and $x,y\in X_e$ are adjacent if and only if $(x,y)\in R_2$.
By Lemma \ref{lem:Ue_module_hypercube} the space $\C^{X_e}$ is a $\U_e$-submodule of $\C^X$.
Define
$$
\T_e={\rm Im}\left(
\U_e\to \End(\C^{X_e})
\right)
$$
where $\U_e\to \End(\C^{X_e})$ is the representation of $\U_e$ into $\End(\C^{X_e})$ corresponding to the $\U_e$-module $\C^{X_e}$.
Note that $\T_e$ is the {\it Terwilliger algebra} of $\frac{1}{2}H(D,2)$ with respect to $\emptyset$ \cite[Theorem 6.4]{halved:2023}.

\begin{lem}
[Corollary 6.6(i), \cite{halved:2023}]
\label{lem:irrUemoduleXe}
The pairwise non-isomorphic irreducible $\U_e$-submodules of $\C^{X_e}$ are
\begin{align*}
&L_{D-2k}^{(0)}
\qquad
\hbox{for all even integers $k$ with $0\leq k\leq \floor{\frac{D}{2}}$};
\\
&L_{D-2k}^{(1)}
\qquad
\hbox{for all odd integers $k$ with $1\leq k\leq \floor{\frac{D-1}{2}}$}.
\end{align*}
\end{lem}

By Lemma \ref{lem:R_module_hypercube} the space $\C^{X_e}$ is an $\Re$-submodule of $\C^X$. Define
$$
\R_e={\rm Im}\left(
\Re\to \End(\C^{X_e})
\right)
=\{M|_{\C^{X_e}}\,|\, M \in \R\}
$$
where $\Re\to \End(\C^{X_e})$ is the representation of $\Re$ into $\End(\C^{X_e})$ corresponding to the $\Re$-module $\C^{X_e}$. 
Both $\T_e$ and $\R_e$ are semisimple subalgebras of $\End(\C^{X_e})$. We end this paper with a comment on $\T_e$ and $\R_e$:

\begin{thm}
\label{thm:Te&Re}
The algebra $\R_e$ is a subalgebra of $\T_e$. Moreover the following conditions are equivalent:
\begin{enumerate}
\item $\R_e=\T_e$.

\item $D$ is odd.
\end{enumerate}
\end{thm}
\begin{proof}
It follows from Proposition \ref{prop:sharp&Ue} that $\R_e$ is a subalgebra of $\T_e$.

(i) $\Rightarrow$ (ii): Suppose on the contrary that the positive integer $D$ is even.
By Lemma \ref{lem:irrUemoduleXe} there is an irreducible $\U_e$-submodule of $\C^{X_e}$ isomorphic to $L_D^{(0)}$. By Theorem \ref{thm:irr_Ln0}(iii), (iv) the $\Re$-module $L_{D}^{(0)}$ is reducible, a contradiction to (i).

(ii) $\Rightarrow$ (i): Suppose that (ii) holds. Applying Theorems \ref{thm:irr_Ln0} and \ref{thm:irr_Ln1} to Lemma \ref{lem:irrUemoduleXe}, as $\Re$-modules,  the pairwise non-isomorphic irreducible $\U_e$-submodules of $\C^{X_e}$ are isomorphic to the irreducible $\Re$-modules
$
R_{\frac{D-1}{2}-k}(-\frac{1}{4},-\frac{1}{4},-\frac{1}{4})$
for all integers $k$ with $0\leq k\leq \frac{D-1}{2}$.
By Proposition \ref{prop:irr_Ln0&Ln1} these $\Re$-modules are mutually non-isomorphic.
Since $\R_e$ is a subalgebra of $\T_e$ and both are semisimple subalgebras of $\End(\C^{X_e})$, the statement (i) follows.
\end{proof}

\subsection*{Funding}
The research is supported by the National Science and Technology Council of Taiwan under the projects NSTC 112-2115-M-008-009-MY2 and MOST 111-2115-M-A49-005-MY2.

\subsection*{Data Availability}
No data was used for the research described in the article.

\section*{Declarations}
\subsection*{Conflict of interest}
The authors have no conflict of interest to declare.

\bibliographystyle{amsplain}
\bibliography{MP}

\end{document}